\documentclass[12pt]{amsart}

\usepackage{amssymb,amsmath,enumerate}
\usepackage[utf8]{inputenc}
\usepackage[T1]{fontenc}
\usepackage[pdftex]{graphicx}

\theoremstyle{definition}
\newtheorem{proposition}{Proposition}
\newtheorem*{corollary}{Corollary}
\newtheorem*{pf}{Construction}
\newtheorem*{remark}{Remark}

\begin{document}

\title{Georg Mohr's ``Euclides Danicus'' -- Preliminary Version}
\author{Ricardo Bianconi \and (translation and commentary)}

\maketitle

\tableofcontents

\listoffigures

\part{Introduction}

\section*{Foreword}

We present here a translation of the Dutch version of Georg Mohr's \emph{Euclides Danicus}, published in Amsterdam in 1672, 125 years before \emph{La Geometria del Compasso} of Lorenzo Mascheroni.

He published two versions of this work, one in Danish, his mother tongue, and another in Dutch, his adopted language. Both are dated January 31st, 1672. There are a few significant differences between the two versions which are pointed out in the text.

The Danish version was translated into German by Julius P\'al and publish together with an article by Johann Hjelmslev in 1928. Mohr claims in the preface to the reader that  ``\emph{\dots I translated this (on request of good Friends) from my mother's language into Dutch \dots}''. In some places he included new constructions which are not present in the Danish version.

\section{Georg Mohr}

There is scarce information about Georg (J\o rg) Mohr, mostly from short references in letters by Tschirnhauss (his friend), Leibniz, Oldenburg, and also some family records. We summarize some information from Johannes Hjelmslev's \emph{Beitr\"age zur Lebensbeschreibung von Georg Mohr (1640-1697)}, \cite{hjelmslev}, who complains the there is very little about his life in Denmark \cite[p. 20]{hjelmslev}. His sources are some few letters and a brief account on Georg Mohr by his son Peter Georg Mohrenthal (quoted in \cite[pp. 14-15]{hjelmslev}.

\begin{enumerate}
\item
He was born in Copenhagen, 1st April 1640, son of the tradesman David Mohrendal,  and died in Kieslingswalde (now, the Polish village of S\l awnikowice) , near Goerlitz, Germany, 26th January 1697.
\item
He left Copenhagen to go to Amsterdam in 1662, probably to study with Huygens. He became friend of Ehrenfried von Tschirnhaus (1651--1708).
\item
 In 1683, when he was already back in Copenhagen, he still had correspondence with Pieter Van Gent (1940--1693) and Ameldonck Bloeck (1651-1702) in Amsterdam\footnote{Handwritten letter from Mohr to Tschirnhaus, Copenhagen on January 30, 1683.}, who had belonged to the circle surrounding Spinoza.
\item
 He fought in the Dutch-French conflict of 1672-3 and was taken prisoner.
\item
He published the \emph{Euclides Danicus} in 1672 in Danish and Dutch.
\item
In 1673 he published the \emph{Compendium Euclidis Curiosi}, where he proves that any point construed with ruler and compass can be construed with a ruler and a compass with a fixed opening. This book was published anonimously but a letter from Oldenburg attributes this work to Mohr.
\item
He married Elizabeth von der Linde of Copenhagen on 19 July 1687.
\item
In 1695 he took a job with Tschirnhaus (in his museum in Kieslingwalde which he kept till his death two years later.
\end{enumerate}

\section{Mohr-Mascheroni Theorem}

Both Mohr and Mascheroni described several geometrical constructions with compass alone and a small subset of those constructions constitute the proof of Mohr-Mascheroni Theorem.

We list the relevant constructions with pointers to the texts. These are some of the propositions of the first part of his book.

The proof of the theorem requires to construct intersections of two circles, a circle and a line (given by two of its points) and two lines (again each given by two of its points). The reader can find a detailed description of the constructions together with a partial translation into French of the First Part of Mohr's book in \cite{escofier-laouenan}.

\begin{enumerate}
\item
If we know the centre of the circle and one of its points then there is no secret on drawing the circle. Otherwise, only one diameter (a pair of points) is given: the centre is the midpoint, construed in Proposition 15. Also, the circle may be given by three of its points: this case is treated in Proposition 33.
\item
The intersection of a circle with a straight line has two cases:
\begin{enumerate}[(a)]
\item
the line goes through the centre of the circle, treated by Proposition 32;
\item
the line does not goes through the centre of the circle, treated in Propositions 19 together with Proposition 22.
\end{enumerate}
\item
The intersection of two lines are treated in Proposition 31 (to find the fourth proportional of three given segments).
\end{enumerate}

\section{On the Book}

The book is written in a terse style in which the figures play a necessary role. Mohr refers to points in the figures without describing them in the text (as, for instance, we never see something like \emph{Let $C$ be the point such that \dots}, although we see things like \emph{resulting in $C$}). Mascheroni's book contains the problems, the solution (construction, or synthesis) and a proof of its correctness (analysis).

Also Mohr states the propositions (problems) and but for a unique exception (Proposition 21 of the First Part) he gives no proof of the correctness of the described construction.

The book is divided in two parts.

The First Part deals with all the fundamental constructions of plane geometry needed to solve more complex problems, and ends with some applications to problems of areas of figures, totalling 54 propositions.

The Second Part contains some more intricate constructions, dealing with intersections (Propositions 1--6), tangency (Propositions 7--10), partitions (Propositions 11--16), centre of gravity (Proposition 17) , Snellius-Pothenot Problem\footnote{Given three points $A$, $B$ and $C$; whose observer stood in $E$ and found the angle $AEB$ equal to the arc $FG$ and the angle $BEC$ equal to the arc $GH$. One asks for the point $E$  where the observer was.} (Proposition 18), perspective (Propositions 19--22) and sundials (Propositions 23 and 24).

There are some small and other significant differences between the Danish and the Dutch texts. A first big difference is a presence of a foreword in the Dutch version (\emph{To the Art Loving Reader}), and absent in the Danish version, in which Mohr explains his motivations and what he intend to show. Most of it is contained in the afterword to the Danish version, reproduced below (translated from the German translation by Julius P\'al).

\begin{center}\textbf{The Afterword of the Danish version} \end{center}

I hope that what is contained in this little work can be reproduced by anyone who has some knowledge of the Elements of Euclid (because it doesn't appear to be difficult), as well as the evidence (and the missing Elements of Euclid), which are omitted here, can easily be provided by oneself. But if there are some who are not pleased by this manner of solution because it is easier to carry out with a ruler and compass, they should know that I am well aware of this; but my only concern was to examine some of the nature of the circle, whether its characteristic was that it could be used to solve the planimetric constructions (without using a ruler) by cutting circles as detailed here. I beg the willing reader that he will accept this with love, especially if a mistake should have occurred, considering that all of our works are imperfect.

\section{The sources of the figures}

The figures which make an essential part of the text are scanned from the original source.

\begin{figure}[!ht]
\begin{center}
\includegraphics[width=\textwidth]{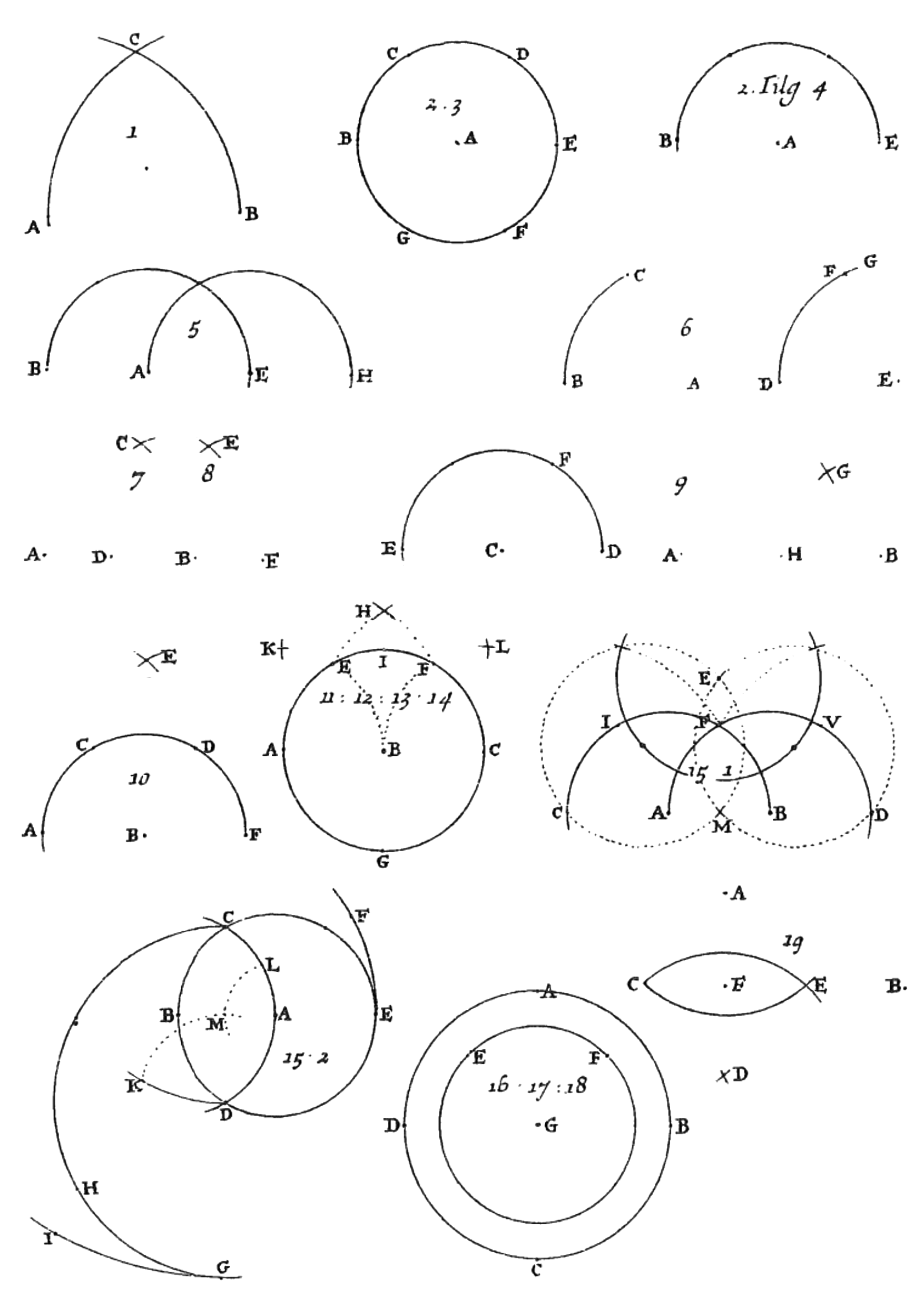}
\end{center}
\caption{Propositions 1 to 19. }\label{fig-prop-1-19}
\end{figure}

Figure \ref{fig-prop-1-19}, page \pageref{fig-prop-1-19}:

By Georg Mohr - Scan of a reprint of the 1672 book, Public Domain, https://commons.wikimedia.org/w/index.php?curid=68212559

\begin{figure}[!ht]
\begin{center}
\includegraphics[width=\textwidth]{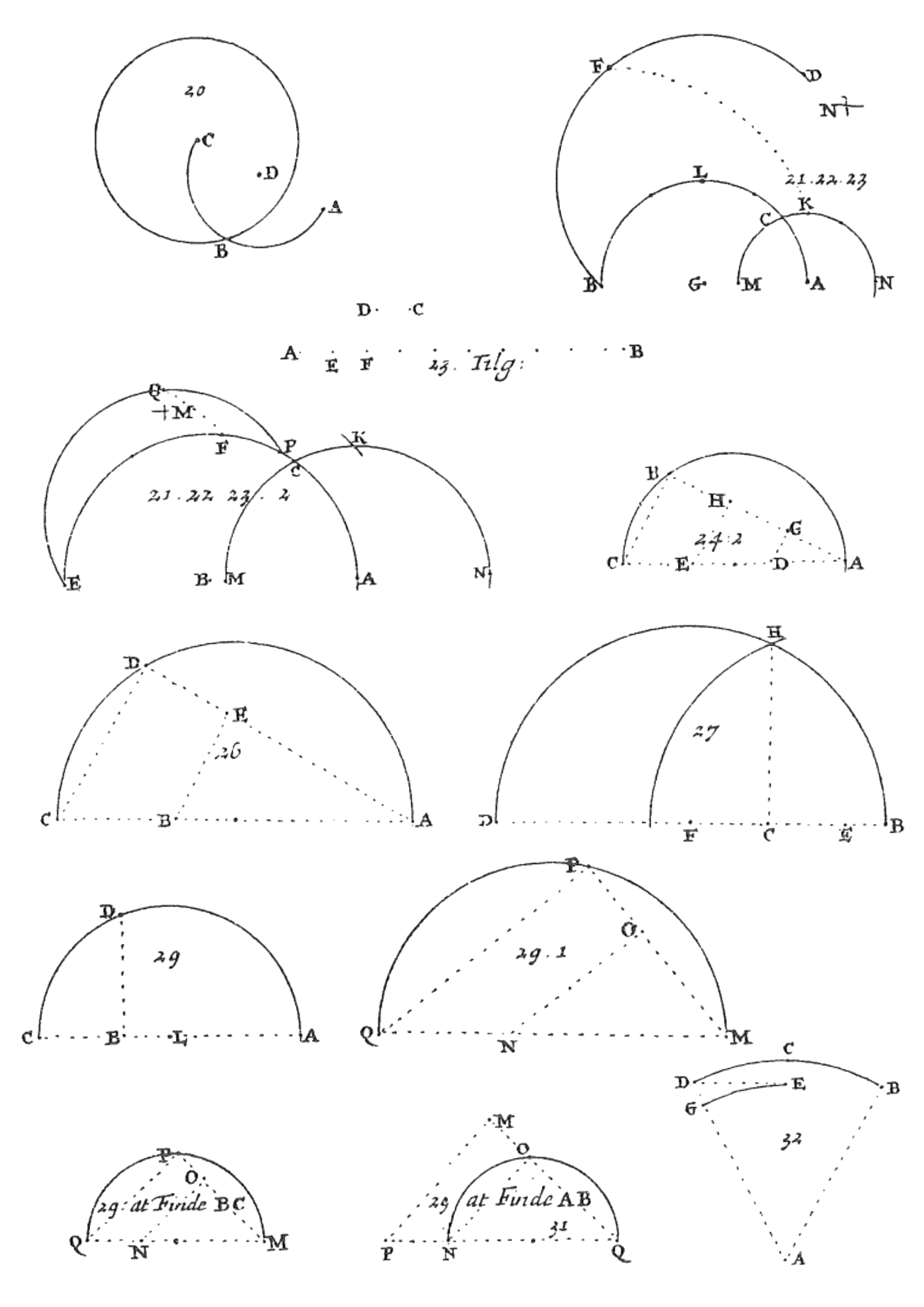}
\end{center}
\caption{Propositions 20 to 32. }\label{fig-prop-20-32}
\end{figure}

Figure \ref{fig-prop-20-32}, page \pageref{fig-prop-20-32}:

By Georg Mohr - Scan of a reprint of the 1672 book, Public Domain, https://commons.wikimedia.org/w/index.php?curid=68212560

Figure \ref{fig-prop-33-45}, page \pageref{fig-prop-33-45}:

By Georg Mohr - Scan of a reprint of the 1672 book, Public Domain, https://commons.wikimedia.org/w/index.php?curid=68212561

Figure \ref{fig-prop-46-54}, page \pageref{fig-prop-46-54}:

By Georg Mohr - Scan of a reprint of the 1672 book, Public Domain, https://commons.wikimedia.org/w/index.php?curid=68212563

Figure \ref{fig2-prop-1-12}, page \pageref{fig2-prop-1-12}:

By Georg Mohr - Scan of a reprint of the 1672 book, Public Domain, https://commons.wikimedia.org/w/index.php?curid=68212565

Figure \ref{fig2-prop-13-18}, page \pageref{fig2-prop-13-18}:

By Georg Mohr - Scan of a reprint of the 1672 book, Public Domain, https://commons.wikimedia.org/w/index.php?curid=68212564

Figure \ref{fig2-prop-19-22}, page \pageref{fig2-prop-19-22}:

By Georg Mohr - Scan of a reprint of the 1672 book, Public Domain, https://commons.wikimedia.org/w/index.php?curid=68212568

Figure \ref{fig2-prop-23-24}, page \pageref{fig2-prop-23-24}:

By Georg Mohr - Scan of a reprint of the 1672 book, Public Domain, https://commons.wikimedia.org/w/index.php?curid=68212569


\part{The Translation}

\setcounter{section}{0}

\section{Frontpages}

\begin{center}
\textbf{Illustrious, Most Powerful King}

\textbf{Most Merciful Sir and Heir-King}
\end{center}

The innate love that a man posses to his Fatherland does not cool down easily, but he himself tries to help with his service, in whatever way, specially in what he understands best. Therefore I could not refuse myself to bring to light, with my inabilities (my duty so demands), this little scientific work, as my first sprout  in our Danish mother language, treating the plane geometry of Euclid and also of various other writers, but solved in a completely different way, which I believe no other author has done in any language. And since the high authority rightly bear the name of Father of the Fatherland, so find me obliged to Your Royal Majesty, my Most Compassionate Heir-King and Sir, to humbly dedicate this Treatise as a token of my very submissive duty. Praying most benignly, Your Royal Majesty please take this my little but sincere labour, and lend the missing lustre to the same. May GOD Almighty grant Your Royal Majesty and Royal family with all desirable happiness, salvation and prosperity, blessing and protection. What wishes,

Your Royal Majesty,

As my Most Compassionate

Heir-King and Sir,

Amsterdam, 31st January, 1672.

Most Compassionate, faithful subject,

Georg Mohr.

\bigskip

\begin{center}\textbf{To the Art Loving Reader}\end{center}

One knows that instruments are nothing but auxiliary tools for the geometric constructions, when at first they are deduced by mathematical reasoning, as one sees in the study of Euclid's plane geometry, where one uses rule and compass to trace straight lines and circles.
So I have undertaken the task of studying the nature and possibilities of the compass alone, making the constructions which I found feasible with intersections of circles, and such that some mathematicians still had doubts about being possible, others even considered impossible, and some required that I brought to light my findings.
Therefore I found myself forced to show my latest discoveries given in this Treatise, divided in two parts. The first one contains in most part Euclid, its constructions from first Elements, and the second one contains various things, although I have not  many examples of each case because I thought them to be superfluous.
As for the exposition of the Elements, they follow one another in such a way that the next one depends on the previous ones, and therefore the second one is Proposition 15 of Euclid's Fourth Book\footnote{To construct a regular hexagon.}, and so on, and similarly as one can see, what distinguishes this work from Euclid's is that what he does with ruler and compass, it uses only the compass, what I leave to the wise reader to judge (and also to those who let their minds fly). 
And although it is possible to solve most of the Elements in different ways, so I have said: I hope that it can be done for each one (if it doesn't seem heavy right now) who have little knowledge of Euclid's Elements; and also their evidence (and the missing Elements of Euclid) that I have left behind, can now be easily found; but there may be someone who could not be pleased by this way of solving, while it can be done more easily by rule and compass, who must know that this is not unknown to me, but my purpose is only to ponder on the nature of the Compass or to show its potential, as it was said before. 
Here is also a comment on the good-natured Reader, that he would like to keep this little work at its best, and not on the decorative side of writing; all the more that I translated this (on request of good Friends) from my mother's language into Dutch, but please accept the good affection for the doing.
And so herein  a few faults may have crept in, please at best indicate such, and remember that our work is only a piecework.

\section[First Part]{First Part: Deals with Euclid's plane constructions, comprehended in the first book of Euclid's [\emph{Elements}]}

\begin{proposition} 
Describe an imaginary equilateral triangle over a given imaginary line $AB$.
\end{proposition}

\begin{pf}
Describe two arcs from $A$ and $B$ with length\footnote{In the Danish version, \emph{aabning} or \emph{{\aa}bning}, opening, \cite{mohr1928}. Mohr uses \emph{length}, \emph{width} and \emph{opening} with the same meaning.} $AB$ which cut each other at $C$. So $ABC$ is the required imaginary triangle.
\end{pf}

\begin{proposition} 
Inscribe a regular hexagon in a circle whose radius is $AB$.
\end{proposition}

\begin{pf}
Trace a circle with $AB$, and from $B$ with same opening mark six times in the circle. Thus the required [figure] is obtained.
\end{pf}

\begin{corollary}
If one traces an arc with the same opening $AB$, and with this  one marks three times from $B$, then the three points $B$, $A$ and $E$ will lie in the same imaginary straight line.
\end{corollary}

\begin{proposition} 
Inscribe an equilateral triangle in a circle whose radius $AB$ is given.
\end{proposition}

\begin{pf}
Firstly make as in the foregoing and then the imaginary lines $BD$, $DF$ and $BF$ will form the required triangle.
\end{pf}

\begin{proposition} 
With the given width $AB$ trace an imaginary straight line twice as long as $AB$
\end{proposition}

\begin{pf}
Mark in an arc of circle with radius $AB$ three times with the same opening from $B$ up to point $E$. Thus	 $B$, $A$ and $E$ will be in a straight line.
\end{pf}

\begin{proposition} 
With the given width $AB$ trace an imaginary straight line from $B$ as long as three (or more) times as long as $AB$.
\end{proposition}

\begin{pf}
Firstly double the line $AB$ (by 4). Now describe from $E$ an arc of circle with the same opening $AB$ and mark three times from $A$ with the same opening resulting in $H$. Then $BH$ with $AB$ will be in an imaginary straight line, and thus one can make more times.
\end{pf}

\begin{proposition} 
Trace an angle $DEF$ equals a given angle $BAC$.
\end{proposition}

\begin{pf}
With the opening $AB$ trace an arc of circle\footnote{An arc centred at $E$ with opening $ED$ equals $AB$.} $DE$; take the width $BC$ and mark from $D$ to $F$. Thus the angle $DEF$ equals $BAC$.
\end{pf}

\begin{proposition} 
Construct a triangle from three given imaginary lines $AB$, $AC$ and $BC$, the sum of two of them larger than the third one.
\end{proposition}

\begin{pf}
Mark over one of the given lengths $AB$ from $A$ an arc with the other given length $AC$, and the same from $B$ with the third length $BC$ , which will cut each other in $C$, as required.
\end{pf}

\begin{proposition} 
Make a triangle $DFE$ congruent\footnote{Literally, \emph{equal and {equiangular}.}} to a given triangle $ABC$.
\end{proposition}

\begin{pf}
Make $DE$ equal to $AC$, and $DF$ equal to $AB$, as well as $EF$ equals $CB$. Then both triangle are congruent to each other.
\end{pf}

\begin{proposition} 
Given two imaginary straight lines $AB$ and $CD$, where $AB$ is twice as large as $CD$, one desires to make a circle over the length $AB$ whose radius is the shortest length $CD$.
\end{proposition}

\begin{proposition} 
From the endpoint $B$ from two given points $A$ and $B$ trace a perpendicular.
\end{proposition}

\begin{pf}
Describe a semicircle with $AB$ (by the Corollary of 2) and make over $CD$ an equilateral triangle (by Proposition 1). Thus $EB$ will be perpendicular to $AB$.
\end{pf}

\begin{proposition} 
Inscribe a square  in a circle whose radius $AB$ is given.
\end{proposition}

\begin{pf}
Make semicircles (by the Corollary of 2) with the lengths $EC$ equal $AF$. With these make over $AC$  an isosceles triangle $AHC$, and with the length $BH$ mark from $A$ points $G$ and $I$ or from $C$ the points $G$ and $I$. Then $AI$ is equal to $IC$, $CG$ and $AG$.
\end{pf}

\begin{proposition} 
Inscribe a dodecagon in a circle whose radius $AB$ is given.
\end{proposition}

\begin{pf}
It is as in the previous proposition, and $IE$ equals to $IF$, equals the required side.
\end{pf}

\begin{proposition} 
Make a square over an imaginary straight line $AB$.
\end{proposition}

\begin{pf}
Firstly find $I$ (as in 11), make now $IK$ equal to $BI$ and $AK$ equal to $AB$. Thus $AKIB$ is the required square.
\end{pf}

\begin{proposition} 
Circumscribe a square on a circle whose radius $AB$ is given.
\end{proposition}

\begin{pf}
Inscribe a square in the circle (as in 11), then make a square over $AB$ and $BC$ or $BI$ and $BG$ (by 13). Thus $KLMN$ is the required square.
\end{pf}

\begin{proposition} 
We desire to find the midpoint\footnote{Mohr gives three constructions and Mascheroni gives other five constructions in \cite[Libro Terzo, Problema 66, pp. 31--34]{mascheroni}, making a total of eight distinct constructions.} $M$ of two given points (or an imaginary straight line) $A$ and $B$, so that $A$, $M$ and $B$ lie in an imaginary straight line.
\end{proposition}

\begin{pf}
Describe two semicircles from $A$ and $B$ with opening $AB$ (by Corrollary to 2), and make an isosceles triangle $CED$ over $CD$ (by 7) such that the upstanding sides $CE$ and $DE$ be equal to $CB$ or $AD$. Now describe a circle over $CE$ and [one over] $ED$ (by 9), that will cut themselves in $M$, which is the required midpoint of $AB$.

\medskip

\textsc{Otherwise.}
Describe two semicircles from $A$ and $B$ with opening $AB$ (by Corollary of 2), which will cut themselves in $C$ and $D$. From $B$ with opening $EB$ make the arc $EF$ and from $E$ mark $EF$ equal $AB$; now take the opening $FA$ and mark from $E$ to $L$. Thus $AL$ will be equal to $AM$, the half of $AB$. In the same way find the half of $CD$ resulting in $DK$ equal to $DM$, and the arcs $AL$ and $DK$ will cut themselves at $M$, which is the midpoint of $AB$.

\begin{remark}
After this manner one can find any required division\footnote{For instance, repeat this construction but using a point $E$ such that $BE$ is three times as large as $AB$, then $BM$ will be one third of $AB$. This procedure is repeated in the first construction of Proposition 24 below.}\label{footnote-prop-15}  of $AB$.
\end{remark}

\textsc{Otherwise.}
Make an isosceles triangle $CPD$, whose upstanding sides $CP$, equal to $DP$, each equal to $EB$, and the base $CD$ equal $EC$. Describe a circle on $CP$ and on $DP$ (by 9) which will cut themselves in $M$, the required midpoint.
\end{pf}

\begin{proposition} 
Describe a circle though two given points $A$ and $C$ whose diameter is $AC$.
\end{proposition}

\begin{pf}
By the foregoing divide $AC$ in the middle in $G$, which is the required midpoint.
\end{pf}

\begin{proposition} 
Inscribe a circle in a given square $ABCD$.
\end{proposition}

\begin{pf}
Divide $AB$ or $AD$ in the middle, in $E$ or $F$ (by 15), and construct a square over $AE$ or $AF$ (by 13). Thus the intersection $G$ is the required midpoint.
\end{pf}

\begin{proposition} 
Circumscribe a circle on a given square.
\end{proposition}

\begin{pf}
Divide $AC$ or $BD$ in the middle in $G$, which is the required midpoint.
\end{pf}

\begin{proposition} 
Drop a perpendicular from a given point $A$ above two given points (or an imaginary straight line) $B$ and $C$ (say, $AF$) so that $B$, $F$ and $C$ lie in an imaginary straight line.
\end{proposition}

\begin{pf}
Construct the triangle $BDC$ equal to the triangle $BAC$ (by 8); take the width $AC$ equal to $DC$ and trace from $A$ and $D$ arcs which cut themselves in $E$. Divide now $EC$ in the middle in $F$ (by 15), the required point, or divide $AD$ in the middle and also arrive in $F$.
\end{pf}

\newpage

\begin{proposition} 
Trace from a point $A$ an imaginary line $AB$ that touches a given circle\footnote{The point $C$ is the centre. The circle is assumed already traced.} $BC$.
\end{proposition}

\begin{pf}
Divide $AC$ in the middle in $D$ (by 15), and describe a circle with $AD$ equal $DC$ which will cut the given circle in $B$, which is the required point.
\end{pf}

\begin{proposition} 
One desires to trace an imaginary straight line\footnote{The smallest length $AC$ is given.} $AC$ perpendicular to a given imaginary straight line $AB$ through one of its extremities.
\end{proposition}

\begin{pf}
Describe over the longest line $AB$ a semicircle (by 16), and also
from $A$ describe an arc with the length $AC$ which will cut the bigger one in $C$. Now from $A$ or $B$ inscribe a square in the circle $AB$ whose side is $AL$, equal to $BL$, (by 11). With this length $BL$, from $L$, describe a semicircle $BFD$ (by the Corollary of 2); mark therein mark $BF$ from $D$ in $F$; take now $BF$ and mark from $B$ in $K$, then $AK$ will be perpendicular to $AB$.

\smallskip

If there is doubt about this manner, then I shall add a proof\footnote{This does not appear in the Danish version, \cite{mohr1928}. It is the application of Pythagoras Theorem to the relevant triangles.} by the precept.

Let $AB$ be equal to $a$, $AC$ equal to $b$, $BG$ equal to $GA$, equal to $a/2$. As $BL$ is equal to $AL$, equal to $a^2/\sqrt{2}$, so $BD$ is equal to $a^2\sqrt{2}$. Now $AC$ is marked in the semicircle $BLA$ and its square $b^2$ subtracted from the square of $AB$, equal to $a^2$, results in the square of $BC$, equals to [the square of ] $BF$, equals to $a^2-b^2$; this square subtracted from the square of $BD$, equals $2a^2$, results in the square of $BF$, equal to the square of $BK$, equal to $a^2+b^2$, which is equal to the square of $AB$ plus the square of $AK$, equal to the [square of] $AC$, as required.

\medskip

\textsc{Otherwise.}\footnote{This construction also does not appear in the Danish version, \cite{mohr1928}.}
Describe a semicircle $ECA$ from $B$ with the width $AB$ (by the Corollary of 2), and also from $A$ with the length $AC$ describe an arc which will cut the circle $ECA$ in $C$. Now describe a semicircle $EQP$ over $EP$ (which is the side of an isosceles triangle in the circle $ECA$), and therein mark from $P$ $PQ$ equal to $AC$, and the length $EQ$ mark from $E$ to $F$. Now take the width $AF$ and mark $K$ from $B$. Thus $AK$ is perpendicular to $AB$ as before.
\end{pf}

\begin{proposition} 
Two imaginary straight lines $AB$ and $AC$ are given, which one desires to sum together, lying in an imaginary straight line.
\end{proposition}

\begin{pf}
Draw $AC$ perpendicular to $AB$ (by the foregoing) in $K$; from $K$ inscribe a square in the circle of $AC$ (by 11). Thus $KN$ is a side of the square and $N$, $A$ and $B$ are in an imaginary straight line.
\end{pf}

\begin{proposition} 
Given two unequal imaginary straight lines $AB$ and $AC$, one desires to subtract\footnote{We use also the phrasal verb ``\emph{pull off}'' in several occasions below to hopefully improve readability: for instance, when one subtracts a segment $AB$ from another $CD$ starting at the point $A$, we write pull $CD$ off $AB$ from $A$. Mohr writes \emph{subtract $CD$ from $AB$ from $A$}} the shortest\footnote{Mohr assumes that $AC$ is the shortest line.} from the longest one.
\end{proposition}

\begin{pf}
Draw $AC$ perpendicular to $AB$ (by 21) in $K$, and from $K$ inscribe a square  in the circle of $AC$ (by 11). Then $KM$ is the side of the square and $MB$ is the remainder.
\end{pf}

\begin{corollary}
Thus now one has several points in between a given imaginary straight line $AB$ which lie in a straight line with $AB$. So take any imaginary line $CD$ as one pleases and pull off $AB$, from $A$ or $B$, say\footnote{Literally: ``I take from $A$''.} from $A$, resulting in $E$; with this length $AE$ make any number of times from $E$ (by 5), and thus all these points will lie in a straight line with $AB$. Notice that between two given points or an imaginary straight line one can find several points, and so I put dotted lines in the following figures to make them more intelligible.
\end{corollary}

\begin{proposition} 
One desires to divide a given imaginary straight line $AB$ in as many parts as required, say three parts, so that all the points lie in a straight line with $AB$.
\end{proposition}

\begin{pf}
One can find each part by the second manner as given\footnote{See footnote \ref{footnote-prop-15}, page \pageref{footnote-prop-15}, above.} in 15, seeing it otherwise, as that $BE$ must be so many times longer than $AB$, as one likes to divide $AB$. If one does the rest as it says there, one obtains one of the parts, one can pull it off $AB$ from $A$ or $B$ (as in 23); the rest is self-evident (as in 5)\footnote{In the Danish version ``as in 4th or 5th''.}.

\medskip

\textsc{Otherwise.}
Take any two points $A$ and $D$, making it as long as many times one desires to divide $AB$, as $AC$, say three times $AD$ (by 5); now describe a semicircle over $AC$ (by 16), and therein mark the length $AB$ from $A$; drop perpendiculars $HE$ and $GD$ over $AB$ (by 19). Thus $AB$ is divided as required.

\medskip

\textsc{Otherwise}\footnote{This construction appears only in the Dutch version. Its terse text is somewhat difficult to understand as it is. Let $PQR$ be the isosceles triangle whose base $PQ$ and upstanding side $PR$ are both equal to three times $AB$ and the other side $QR$ equal to twice $AB$. Drop a perpendicular from $R$ to $PQ$, resulting in $S$. The base $PQ$ is divided in two parts $PS$ and $SQ$, the latter being the smaller part, measuring twice the third part of $AB$. Check this by similarity of triangles.}.
Describe an isosceles triangle whose base and one upstanding side as long as many times one desires to divide $AB$, and the other upstanding side equal to twice the length of $AB$. Draw over this imaginary side a circle (by 9) as well on its opposite equal and uniform triangle (by 8 and 19) cut the ground in two parts, the smallest one is always equal to twice the required part; the rest work as in 23 and 5. Thus one obtains what was required.
\end{pf}

\begin{proposition} 
Cut off a required part of a given imaginary straight line $AB$, say one third.
\end{proposition}

\begin{pf}
It is equal to the second manner of the foregoing, resulting in $AG$, which is one third of $AB$.
\end{pf}

\begin{proposition} 
One wants to divide a given imaginary straight line $AD$ in a given ratio as $AB$ to $AC$.
\end{proposition}

\begin{pf}
Describe a semicircle  over $AC$ (by 16); mark therein $AD$ from $A$; drop a perpendicular from $B$ do $AD$ (by 19), resulting in $E$. Thus $AD$ is divided in the required ratio.
\end{pf}

\begin{remark}
If $AD$ is larger than $AB$ plus $BC$, then one can always make them as many times larger than $AD$ by 5.
\end{remark}

\begin{proposition} 
Construct a square equal to a given rectangle $DC$, $CB$.
\end{proposition}

\begin{pf}
Add $DC$ and $CB$ together, as $DB$ (by 22); describe over $DB$ a semicircle (by 16) whose centre is $F$; make $CE$ equal to $CF$ (by 4); take now the width $BF$, equal to $FD$, and with this opening from $E$ cut the circle in $H$. Thus $HC$ is the side of the required square.
\end{pf}

\begin{proposition} 
Find the geometric mean between two given imaginary straight lines $DC$ and $CB$.
\end{proposition}

\begin{pf}
It is equal to the foregoing.
\end{pf}

\begin{proposition} 
Make a rectangle from a square whose side $BD$ is given, and one of the sides\footnote{The biggest side.} of the rectangle is $AB$; one wants to find the other side $BC$.
\end{proposition}

\begin{pf}
Make over $MQ$, equal twice $AB$,  a circle (by 16); pull from $M$, $MN$ equal to $AD$; now from $M$ mark $MP$ (equal to $MN$, equal to $AD$) in the circle; drop on it a perpendicular from $N$ falling on $O$ (by 19); thus $MO$ is equal to $AL$, the radius\footnote{Lit. \emph{Radius} in Dutch, \emph{halffmaaleren} in Danish. Similarity of triangles implies that $MO$ is half of $BC$. Therefore the need of the last part of this construction.}; work the rest (by 23 and 4), resulting in $BC$.

\medskip

\textsc{Otherwise on Finding} $BC$.
Make on $QN$, equal to $AB$, a circle (by 16); pull from $M$, $MN$ equal $BD$ (by 23), resulting in $NQ$, and mark from $M$ in the circle $MP$ equal to $MN$; now drop from $N$ a perpendicular falling on $MP$ at $O$ (by 19); then $MO$ is equal the required side $BC$.

\medskip

\textsc{Otherwise on Finding} $AB$ \textsc{given} $BC$.
Draw $QN$ equal to twice $BC$, and $QP$ equal to twice $BD$ together in a straight line (by 22 or 23); describe on $QN$ a circle (by 16); mark therein $QO$ equal to $BD$; drop now from $P$ a perpendicular on $QO$, resulting in $M$. Thus $QM$ is the required side $AB$. Note that if $QO$ is longer than $QN$, then one must make $QN$ larger that $QP$ by the fifth [Proposition].

\textsc{Otherwise}\footnote{This strange construction is not present in the Danish version and there is no drawing for this construction. The translator could not make out what it really is.}
Describe a semicircle over twice the length $BD$ (by the Corollary to 2); mark therein twice $BC$ from one endpoint of the diameter; make hereupon a triangle whose one upstanding side is equal to the side of the equilateral triangle inscribed in the same circle, and the other upstanding side equal to twice the second side of the rectangle in the same circle. Describe a semicircle on this length (by 9), and also on its opposite  the congruent triangle  (by 8 and 19) which will cut the base in a point. Thus the length from this point to the last named side of the rectangle equals $BA$.
\end{pf}

\begin{proposition} 
Find the third proportional to two given imaginary straight lines.
\end{proposition}

\begin{pf}
It is equal to the foregoing, that is, \emph{as $AB$, equal to $MQ$, is to $BD$, equal to $MN$, as $BD$, equal to $MP$, is to $MO$}.
\end{pf}

\begin{proposition} 
Find the fourth proportional to three given imaginary straight lines.
\end{proposition}

\begin{pf}
Make firstly the rectangle equal to a square (by 27); then make the square equal to a rectangle (by 29); thus there is the required outcome.

\medskip

\textsc{Otherwise}.
Similarly, \emph{as $QN$ is to $QP$, so $QO$ is to $QM$}; but if $QO$ is longer than $QN$, then one can make $QN$ as longer as one wants, so that $QN$ becomes longer than $QO$, (by 5).
\end{pf}

\begin{proposition} 
Divide in the middle in $C$ a given angle or sector $BAD$
\end{proposition}

\begin{pf}
Divide $BD$ in the middle (by 15) in $E$; pull $AE$ from $AD$ (by 23); take $GD$, equal to $EC$, and erect a perpendicular on $BE$ from $E$ (by 21); thus one obtains the required [point].

\medskip

\textsc{Otherwise}.
Find the third proportional to $AD$ plus $AE$ and $DE$, resulting in $EC$, like the foregoing.
\end{pf}

\begin{proposition} 
Describe a circle through three given points $A$, $B$ and $C$.
\end{proposition}

\begin{pf}
Find the fourth proportional to three lines, that is, \emph{as twice $BD$ is to $AB$, so $BC$ is to the radius $AE$, equal to $EC$, of the required circle}.
\end{pf}

\begin{proposition} 
Given an arc of circle $AB$, describe the whole circle.
\end{proposition}

\begin{pf}
Pick any point $C$ in the arc and describe the circle through the three points $A$, $B$ and $C$ by the previous [proposition].
\end{pf}

\begin{proposition} 
Cut an arc of a circle $AB$ in the middle in $C$.
\end{proposition}

\begin{pf}
Make firstly as in the foregoing and then make as in 32, thus resulting in $C$ as required.
\end{pf}

\begin{proposition} 
Inscribe in a given circle whose radius is $ED$ a triangle similar\footnote{Lit. \emph{gelijckhoeckig}, or \emph{equal angles}.} to a given imaginary triangle $ABC$.
\end{proposition}

\begin{pf}
Circumscribe a circle on the three points $A$, $B$ and $C$ (by 33) and trace this circle it in the given one. Now pull $BD$ from $ED$ (by 23), resulting in $BE$, and the same with the remaining two, resulting in $F$ and $G$. Thus the triangle $EFG$ is similar\footnote{lit. \emph{gelijck}.} to the triangle $ABC$.

\medskip

\textsc{Otherwise}.
Draw the shortest side $BC$ from $C$ and $A$ from $CA$, $AC$, and $AB$, resulting in $AR$, $TC$ and $SB$ (as in 23); set from $E$ the perpendiculars\footnote{On opposite sides of $DE$.} to $DE$, $EM$ and $EN$ (as in 10); make on $EN$  and $EM$, from $E$, the sector $PEK$ equal to the sector $RCB$, and the sector $OEV$ equal to the sector $TAS$ (by 6). Drop from $K$ a perpendicular falling on $DE$ at $F$, as well as one from $V$ falling on $DE$ at $G$ (by 19); find now the fourth proportional to three lines (by 31 and 23) resulting in $EQ$ and $EY$ (that is, \emph{as $EK$ is to $EF$, so $ED$ is to $EQ$}, and equally, \emph{as $EV$ is to $EG$, so $ED$ is to $EY$}); double each one (by 4), resulting in $X$ and $W$; thus the triangle $XEW$ is similar\footnote{Lit. \emph{gelijckhoeckigh}.} to the triangle $CBA$.
\end{pf}

\begin{proposition} 
Circumscribe on a given circle (whose radius is $HM$) a triangle similar to a given triangle $ABC$.
\end{proposition}

\begin{figure}[!ht]
\begin{center}
\includegraphics[width=\textwidth]{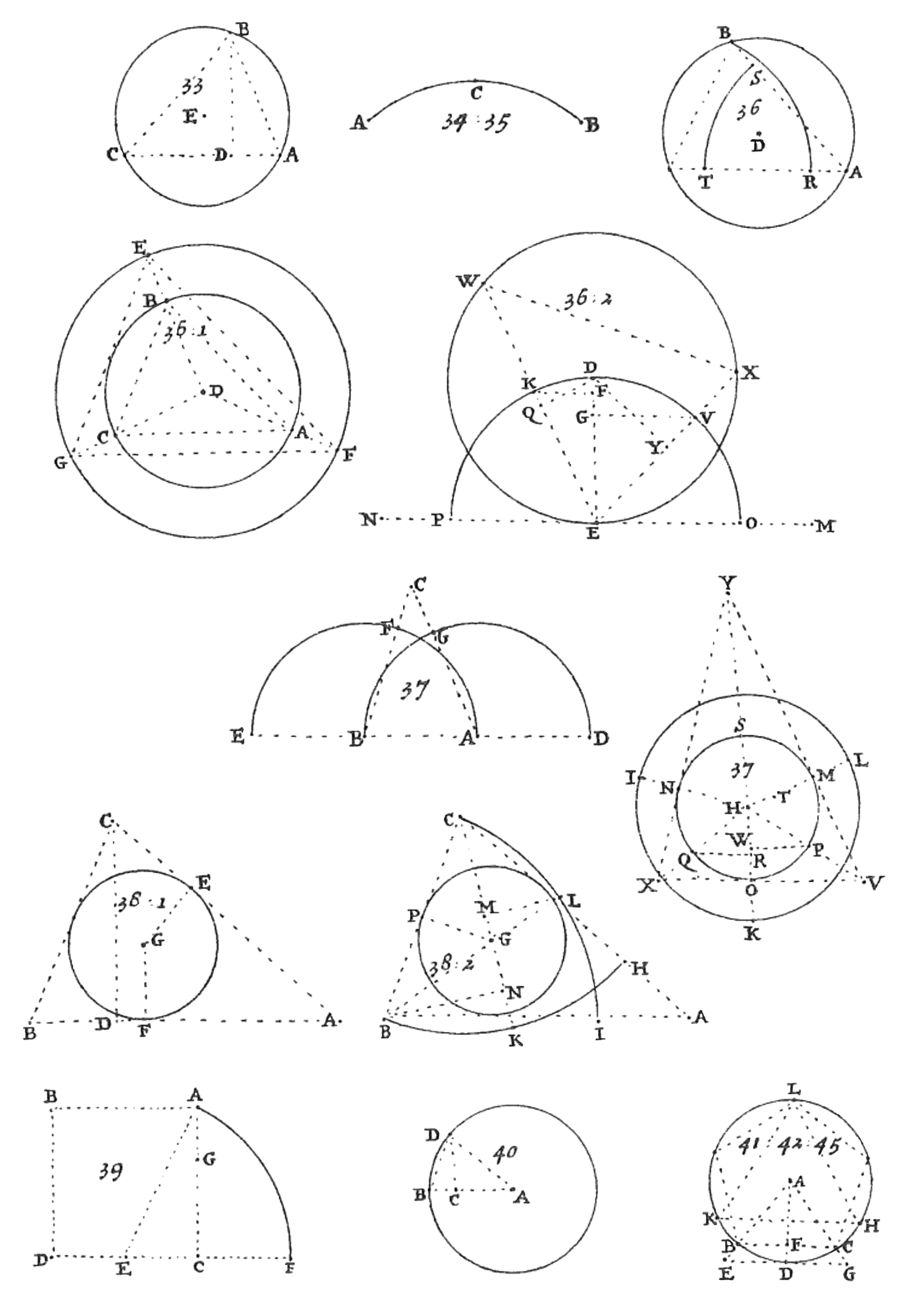}
\end{center}
\caption{Propositions 33 to 45. }\label{fig-prop-33-45}
\end{figure}

\begin{figure}[!ht]
\begin{center}
\includegraphics[width=\textwidth]{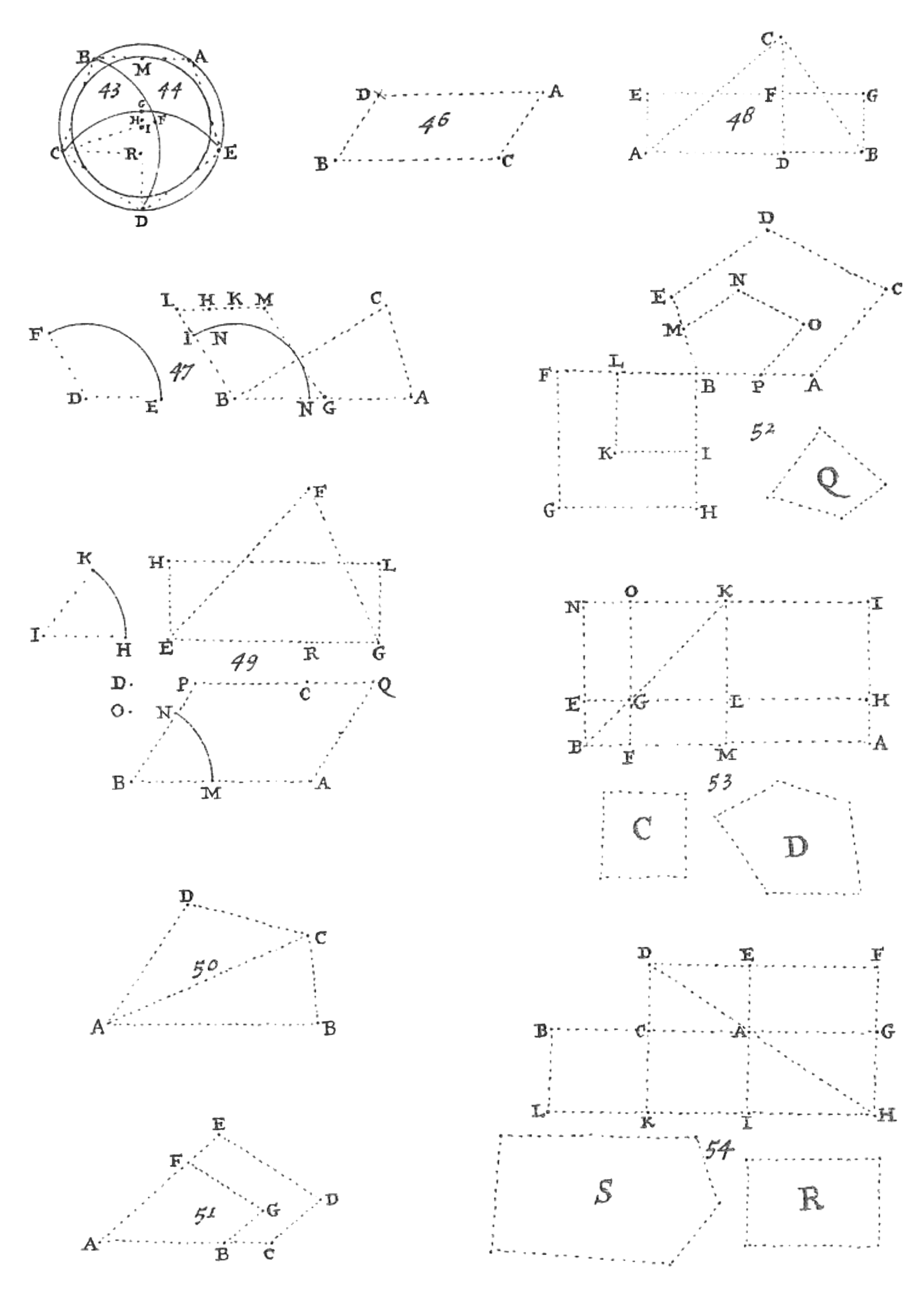}
\end{center}
\caption{Propositions 46 to 54. }\label{fig-prop-46-54}
\end{figure}

\begin{pf}
Pull the shortest side $AB$ off $AC$ and $BC$ from $A$ and $B$, resulting in $GC$ and $FC$; describe from $A$ and from $B$ two semicircles (by the Corollary of 2), $DGB$ and $EFA$, as well as a circle $KLI$ with the width $AB$ and centre in $H$. Make now the sector $KHI$ equal the sector $EBF$ and the sector $KHL$ equal the sector $DAG$ (by 6); now pull $HM$ of $LH$, and $HO$ of $HK$, as well as $HN$ of $HI$, resulting in $LM$, $KO$ and $NI$ (by 23); now divide the sectors $MHO$, $NHO$ and $MHN$ in two equal parts (by 32), resulting in $P$, $Q$ and $S$; from these points drop perpendiculars falling\footnote{Here Mohr is not very precise: $PW,QR\perp HK$ and $ST\perp HL$. See Figure \ref{fig-prop-33-45}, page \pageref{fig-prop-33-45}. He uses $W$ to find the vertex $V$, via similarity of triangles $HWP$ and $HOV$, finding $OV=OM$ as the fourth proportional described above. The same way applies to triangles $HRQ$ and $ HOX$ to find $X$, and $HTS$ and $HMY$ to find $Y$.} on $HK$ and $HL$ at $W$, $R$ and $T$; find now the fourth proportional to three lines, that is, \emph{as $HW$ is to $WP$, so $HO$ is to $OV$, equal to $VM$}; also find $OX$ equal to $XN$ and $MY$ equal to $YN$. Thus the triangle $VYX$ is similar to the given triangle $ABC$.
\end{pf}

\begin{proposition} 
Inscribe a circle in a given triangle $ABC$.
\end{proposition}

\begin{pf}
Drop a perpendicular falling from one angle, as from $C$ over $AB$ at $D$ (by 19); find the fourth proportional to three lines (by 31), resulting in the radius $AG$ (that is, \emph{as the half of the sum of $AB$, $BC$ and $AC$ is to the half of $AB$, so $CD$ is to the radius $GE$}). Now sum two of the lines together, as $AB$ and $AC$ and pull the third $BC$ (by 22 and 23), half the remaining and pull from $A$ off $AB$ or $AC$, resulting in $EC$ or $FB$; now set the radius found earlier $GE$ perpendicular to $AE$ or $AF$, from $E$ or $F$, resulting in $G$, which is the required centre.

\medskip

\textsc{Otherwise}.
Pull the shortest side $BC$ from $B$ and $C$ off $AB$ and
$AC$, resulting in $AI$ and $AH$; divide each of the sectors $IBC$ and $HCB$ in two equal parts at $L$ and $K$ (by 32); now drop perpendiculars from $B$ and $L$ to $KC$, at $N$ and $M$ (by 19); now find the fourth proportional to three lines (by 31), that is, \emph{as the sum of $NB$ and $LM$ is to $MN$, so $NB$ is to $NG$}; now add this\footnote{$NG$.} to $NK$ (by 22); thus $G$ is the centre; now drop a perpendicular from $G$ to one of the sides, for example $BC$, at $P$; then $GP$ is the radius of the the required circle.
\end{pf}

\begin{proposition} 
One desires to divide\footnote{This is the golden ratio.}  a given imaginary straight line $AC$ in $G$, so that the rectangle ($CA$, $GA$) of the divided line and the divided part equals the square of the other part $GC$.
\end{proposition}

\begin{pf}
Draw a square over $AC$ (by 11); now divide $CD$ in the middle in $E$ (by 15); take $AE$ and draw (from $E$) in a straight line with $ED$ (by 22); draw from $C$ a perpendicular $GC$ to $FC$, equal $FC$; thus $AC$ is divided in $G$ as required.
\end{pf}

\begin{proposition} 
One requires to make an isosceles triangle $BAD$, such that the angles at the base are each the double of the third [angle] $BAD$.
\end{proposition}

\begin{pf}
Take any $AB$ as one wants and divide it in $C$ as in the previous [proposition]; make $BD$ equal $AC$ in the circle\footnote{The circle centred at $A$ and with radius $AB$.}; then $DAB$ is the required triangle in the circle $AB$.
\end{pf}

\begin{proposition} 
Inscribe a regular pentagon\footnote{Lit. \emph{gelijcksijdige vijfhoeck}.} in a given circle.
\end{proposition}

\begin{pf}
Make firstly a triangle as in the foregoing; inscribe in the given circle a triangle similar to the one construed (by 36); then $BC$ is the side of the required pentagon.
\end{pf}

\begin{proposition} 
Circumscribe a regular pentagon on a circle.
\end{proposition}

\begin{pf}
Inscribe in the circle a regular pentagon as in the foregoing, whose side in $BC$; divide the sector $BAC$ in the middle (by 32) in $D$, and divide $BC$ in the middle in $F$ (by 15); find now the fourth proportional to three lines (by 31), that is, \emph{as $AF$ is to $BF$, so $AD$ is to $ED$}; double this (by 21 and 4), resulting in $EG$, one of the required sides.
\end{pf}

\begin{proposition} 
Inscribe a circle in a regular pentagon.
\end{proposition}

\begin{pf}
By 32, divide each of the sectors $DCB$ and $CDE$ in the middle, in $E$ and $F$; drop a perpendicular from $F$ and $C$ falling on $GD$ at $H$ and $R$ (by 19); find now the fourth proportional to three lines, that is \emph{as $FH$ plus $RC$ is to $HR$, so $RC$ is to $RI$}; sum this to $RD$ (by 22); then $I$ is the centre , so now one divides $AB$ in the middle in $M$ (by 15), and with the width $IM$ draw the circle, as required.
\end{pf}

\begin{proposition} 
Circumscribe a circle on a regular pentagon.
\end{proposition}

\begin{pf}
Find the centre $I$ (by 43) and with the width $AI$ draw a circle; then one obtains the required one.
\end{pf}

\begin{proposition} 
 Inscribe a regular pentadecagon in a given circle.
\end{proposition}

\begin{pf}
Inscribe an equilateral triangle (by 3) and a regular pentagon (by 41) so that the triangle and the pentagon star in a point $L$; then
$BK$ is equal to $CH$, one side of the pentadecagon, as required.
\end{pf}

\begin{proposition} 
One wants to draw an imaginary straight line from a given point $A$ parallel to another given imaginary straight line $BC$.
\end{proposition}

\begin{pf}
Make an arc from $B$ with length $AC$ and an arc from $A$ with length $BC$, which will cut the first one at $D$. Then $AD$ is parallel to $BC$. Notice that when one wants to draw a given imaginary straight line from $A$ parallel to $BC$, one can do such construction (by 22 or 23).
\end{pf}

\begin{proposition} 
Make a parallelogram $(GBLM)$ equal to a given triangle $ABC$, having one angle equal to a given angle $EDF$.
\end{proposition}

\begin{pf}
From $C$ make $CH$ parallel to $AB$ (by 46), and on
$AB$, from $B$, an angle $NBI$ equal to the angle $EDF$ (by 8); now drop from $B$ a perpendicular to $CH$ (as in 19), resulting in $K$, and the same from $I$ on $KB$, in $N$; find now the fourth proportional to three lines, that is, \emph{as $BN$ is to $NI$, so $BK$ is to $KL$}; sum this to $KC$ from $K$ (by 22); now divide $AB$ in the middle in $G$ (by 15), and from $G$ make a parallel to $BL$, resulting in $M$; then the parallelogram $GBLM$ is equal to the triangle $ABC$, as required.
\end{pf}

\begin{proposition} 
Make a rectangle $(BAEG)$ equal to a given triangle $ABC$.
\end{proposition}

\begin{pf}
Drop from $C$ a perpendicular to $AB$ in $D$ (by 19), and divide $CD$ in two equal parts in $F$ (by 15); now draw from $F$ a parallel with $AD$ and $DB$ (by 46); then $BAEG$ is the required rectangle.
\end{pf}

\begin{proposition} 
Make a parallelogram $(BPQA)$ over a given imaginary line $AB$, equal to a given triangle $EFG$, having an angle equal to a given angle $HIK$.
\end{proposition}

\begin{pf}
Make a rectangle $EGLH$ equal to the triangle $EFG$ by the foregoing; find now the fourth proportional to three lines, that is, \emph{as $AB$ is to $EG$, so $GL$ is to $BD$}; erect with this a perpendicular to $AB$ from $B$ (by 21); now make from $D$ a parallel with $AB$, resulting in $C$ (by 46); make the angle $MBN$ from $B$ over $AB$, equal to the angle $HIK$, and from $N$ drop a perpendicular falling on $BD$ at $O$ (by 19); find now the fourth proportional to three lines, that is, \emph{as $BO$ is to $ON$, so $BD$ is to $DP$}; pull this off $DC$ from $D$ (by 23), resulting in $CP$; 
now make a parallel with $BA$ from $P$, resulting in $PQ$; then $BPQA$ is the required parallelogram.
\end{pf}

\begin{proposition} 
Make a parallelogram equal to a rectilinear figure having an angle equal to a given angle.
\end{proposition}

\begin{pf}
Make a square equal each of the two triangles $ADC$ and $ABC$ (by 48 and 27), then modify these two squares into one square\footnote{Here Mohr is using Pythagoras Theorem disguised as Proposition 21.} (by 21), and do the rest as in 47 or 49.
\end{pf}

\begin{proposition} 
Make a rectilinear figure $ABGF$ over a given straight line $AB$, similar to a given rectilinear figure $ACDE$. 
\end{proposition}

\begin{pf}
Find the fourth proportional to three lines, that is, \emph{as $AC$ is to $AD$, so $AB$ is to $AG$}; pull this off $AD$ from $A$ (by 23), resulting in $DG$; find again the fourth proportional to three lines, that is, \emph{as $AD$ is to $AE$, so $AG$ is to $AF$}; pull this off $AE$ from $A$, resulting in $EF$; thus the required is obtained.
\end{pf}

\begin{proposition} 
Describe a rectilinear figure $(BMNOP)$ equal to a given figure $Q$ and similar to another given figure $BEDCA$.
\end{proposition}

\begin{pf}
Make each of the two rectilinear figures into a square (by 48 and 27); mark their sides $BF$ and $BL$ from $B$ in a straight line with $AB$ (by 22); find now the fourth proportional to three lines (by 31), that is: \emph{as $BF$ is to $BE$, so $BL$ is to $BM$}; now pull this off $BE$ from $B$ (as in 23), resulting in $ME$; one equally finds $ND$, $OC$ and $PA$, as required.
\end{pf}

\begin{proposition} 
On a given straight line $BA$ one wants to make two parallelograms so that one $(BEFG)$ is similar to a given parallelogram $C$ and the other $(AFGH)$ is as big as a given rectilinear figure $D$, but the given rectilinear figure is not bigger than that parallelogram, which can be made on the half of the line [$AB$] and similar to the given parallelogram.
\end{proposition}

\begin{pf}\footnote{The German translation of this construction is not very faithful to the Danish original, which is close to the Dutch version.}
Divide $AB$ in the middle at $M$ (by 15); make a parallelogram $MBNK$ on $BM$ similar to $C$ (by 51); then make the parallelogram $LGOK$ equal to a rectilinear figure (which is the difference of the parallelogram $MBNK$ and the given rectilinear figure $D$) and similar to the given parallelogram $C$ (by 52); now from $M$ make a parallel with $LG$ and likewise from $L$ a parallel with $MA$ (by 46); then the parallelogram $AFGH$ is equal to the given rectilinear figure $D$, and the parallelogram $FBEG$ is similar to the parallelogram $C$, as required.
\end{pf}

\begin{proposition} 
Build upon an imaginary straight line $BA$ a parallelogram $(LBGH)$ equal to a given rectilinear figure $S$, and outside another one similar to a given parallelogram $R$.
\end{proposition}

\begin{pf}
Divide $BA$ in the middle at $C$ (by 15) and make a parallelogram $ACDE$ on $AC$, similar to $R$ (by 51); now make a parallelogram $HKDF$ as large as the parallelogram $ACDE$ together with the given rectilinear figure $S$ (by 52), and on the half line $CB$ make a parallelogram $KLBC$ equal to $IKCA$; now draw from $H$ a parallel with $AI$ (by 46); then the parallelogram $HLBG$ is as large as the given rectilinear figure $S$ and the parallelogram $HIAG$ falling outside is at the same angle as the given parallelogram $R$, as required.
\end{pf}


\newpage

\section{The Figures for the First Part}

\clearpage

\section[Second Part]{Second Part: Gives and introduction. to various constructions, as intersections, tangents, partitions, perspective and sundials.}

\setcounter{proposition}{0}

\begin{proposition} 
Given the point $C$ outside the semicircle\footnote{Over $AB$.}, extending from $AB$ in an imaginary straight line (but such that $AC$ is shorter than the radius); now one wants to draw a straight line from $C$, which intersects the circle in $E$ and ends on the circle  in $F$ in such a way that $EF$ should be twice as long as $AC$.
\end{proposition}

\begin{pf}
Find the geometric mean between the sum of the diameter with $AC$ and $AC$ (by 28 of the first part); draw from the resulting $AC$ (by 23, First Part); draw the rest from $C$ to $E$. Now transfer the centre $G$ on the other side to $H$; draw an arc from $H$ with $EH$ equal to $EG$ that intersects the circle in $F$ (by 19, First Part); then $EF$ is twice as large as $AC$ and lies on an imaginary straight line with $CE$.

\medskip

\textsc{Otherwise}\footnote{This second construction appears only in the Dutch version.}
Draw the double of $AC$ anywhere in the semicircle; drop on it a perpendicular from the centre $G$ (by 4 and 19, First Part); with this length draw a circle from $G$ and trace the imaginary  line from $C$ and tangent to this circle (by 20 and 4, First Part); the required line.

If one wants to find $CF$ three times as long as $CE$ from $C$, one proceeds as follows: find the geometric mean between the third part of the diameter plus the third part of $AC$, and $AC$ (by 24 and 28, First Part), resulting in $CE$; draw this from $C$ three times resulting in $F$; then $CF$ is the required imaginary straight line.
\end{pf}

\begin{proposition} 
From a given point $C$, outside the semicircle $AFB$, to draw an imaginary straight line $CF$ which cuts through the circle in $E$ and ends on the circle in $F$ such that the rectangle $CE$, $EF$ is as large as the rectangle $CA$, and $AG$, equal to the radius of the circle.
\end{proposition}

\begin{pf}
Find the geometric mean $CE$ between $GC$ and $AC$, after which one can find $EF$ as in the foregoing and the requirement is met.
\end{pf}

\begin{proposition} 
Given an equilateral triangle $CAB$ with the parallel $AP$ from $A$ to $BC$; one wishes to draw a straight line from $B$, which intersects $AC$ in $F$ and ends on the parallel $AP$ in $E$, so that $EF$ should be equal to a given imaginary straight line $S$.
\end{proposition}

\begin{pf}
Divide $BC$ in the middle at G (by 15, First Part); draw $GL$ equal to $S$ from $G$ perpendicular to $BG$ (by 21, First Part); pull off $BG$ from $B$ resulting in $ML$; now pull $ML/2$ off the geometric mean between the sum of $ML/4$ with $BC$ and $ML$, and pull the resultant from $A$ off $AC$, resulting in $F$; now add $S$ to $FB$ from $F$ (by 22, First Part), ending up in $E$, as required.
\end{pf}

\begin{proposition} 
Given a square $ABCD$, of which $AD$ is the extension of $DL$ along an imaginary line; one wants to draw an imaginary straight line from $B$, which intersects $DC$ in $K$ and ends in $I$ on the extension $AL$ in such a way that $KI$ should be equal to a given imaginary straight line $R$.
\end{proposition}

\begin{pf}
From A,  set $AE$, equal to $R$, perpendicular to AD (by 21, First Part)\footnote{The Dutch version points wrongly to Proposition 22 here. The Danish version points correctly.}; pull $AD$ from $D$ off $ED$ (by 23, First Part), resulting in $EF$; then pull $EF/2$ off the geometric mean between the sum of $EF/4$ with $AD$, and $EF$; set the result  perpendicularly to $AD$ from $D$, resulting in $K$; now one adds the given line $R$ to $KB$ from $K$, ending up in the point $I$, as required.
\end{pf}

\begin{proposition} 
The point $A$ is given on the given imaginary straight line $EV$; one wants to find a point $O$ between $A$ and $E$, so that the square from $AO$ is to the rectangle determined\footnote{Lit. \emph{begrepen}, past participle of \emph{begrijpen}, ``comprehend''.} by $OE$ and a given $AV$, as $R$ is to $S$, or $AI$ to $AV$.
\end{proposition}

\begin{pf}
Lift from $A$, $AY$ equal to $AE$, perpendicular to $AV$, as well as  from $I$, $RI$ equal to $AI$, perpendicular to $VI$ or $AI$ (by 21, First Part); describe a circle on $YR$ (by 16, First Part), whose centre is $S$; drop from $S$ a perpendicular to $IA$ (by 19, First Part), falling in $T$; now pull $ST$ off $SW$, equal to $SR$, obtaining $TW$ and $TX$ (by 22 and 23, First Part); now find the geometric mean between $TX$ and $TW$ as $TO$ (by 28, First Part); then the point $O$ is on $AE$, as required.

\medskip

\textsc{Otherwise}.
Find the geometric mean $AK$ between $EA$ and $AI$; now halve $IA$ in $M$ (by 15, First Part), pull $MO$ equal to $MK$ from $M$ off $ME$ (by 23, First Part), the result is $O$, the desired point.
\end{pf}

\begin{proposition} 
Given a circle $ADB$, one wishes to find a point $C$ on the circle, so that if one draws an imaginary straight line from this point, which intersects the diameter [$AB$] (in $E$) and ends on the circle in $D$, then $AD$ is to $DB$, as $AE$ to $EB$.
\end{proposition}

\begin{pf}
From A make a square (by 11,  First Part) whose side $AC$ is equal to $BC$, then C is the desired point; now choose an arbitrary point E on the diameter AB (by 23, First Part); find the fourth proportional to three lines, that is: \emph{as $CE$ is to $AE$, so $EB$ is to $ED$} (by 31 and 22, First Part); now draw the imaginary straight lines from $A$ and $B$ to $D$, then these are as required.
\end{pf}

\begin{proposition} 
Given two circles, whose radius are $AB$ and $CD$; one wishes to find a point $F$ (which should lie with $AB$ and $CD$ on an imaginary straight line), so that if one draws an imaginary straight line through it, it should touch the two circles (in $G$ and $H$).
\end{proposition}

\begin{pf}
Find the fourth proportional to three lines (by 31, First Part), that is, \emph{as $CD$ subtracted from $AB$ is to $BE$, so $CD$ is to $EF$}; now add $EF$ from $E$ to $AE$ (by 22, First Part), resulting in $F$, the desired point; so now one draws an imaginary straight line from this point, which touches the two given circles (by 20, First Part), arising $G$ and $H$, as required.
But\footnote{This last part of the construction appears only in the Dutch version.} if one requires that the point $F$ lie in between both circles $AB$ and $CD$, then one finds it thus: similarly, \emph{as $AB$ plus $DC$ is to $BC$, so $CD$ is to $CF$}, or \emph{as $AB$ plus $DC$ is to $BC$, so $AB$ is to $BF$} (by 31 and 22, First Part), and the one works the rest as before.
\end{pf}

\begin{proposition} 
In a given square $ABCD$, a quarter circle $BDA$ is drawn from $B$ with the width $BD$; one wishes to describe the two (largest) circles, such as $QM$ and $HG$, so that they should touch the square (in $M$ and $N$; and in $G$ and $P$) and also the arc (in $E$).
\end{proposition}

\begin{figure}[!ht]
\begin{center}
\includegraphics[width=\textwidth]{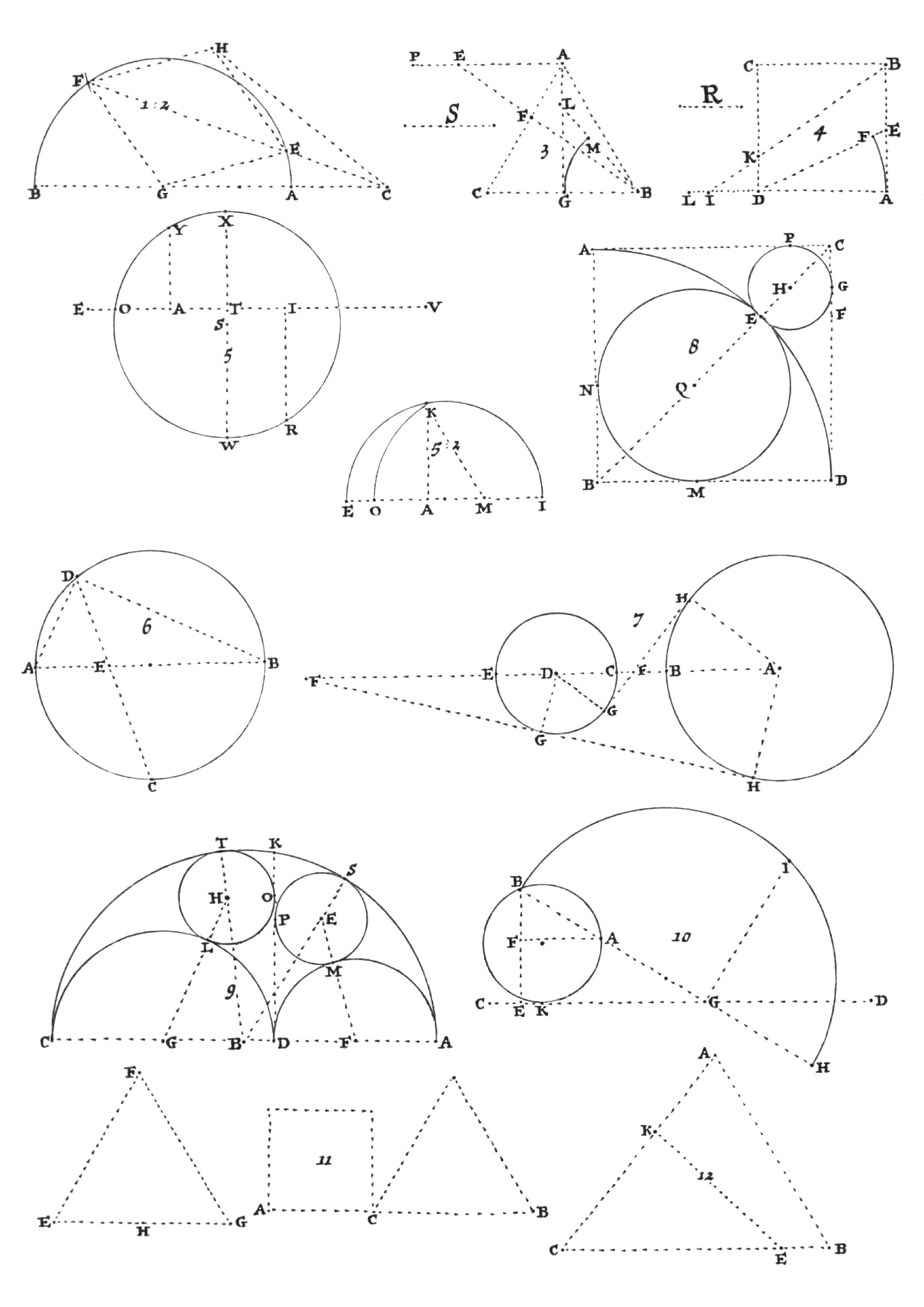}
\end{center}
\caption{Propositions 1 to 12. }\label{fig2-prop-1-12}
\end{figure}

\begin{figure}[!ht]
\begin{center}
\includegraphics[width=\textwidth]{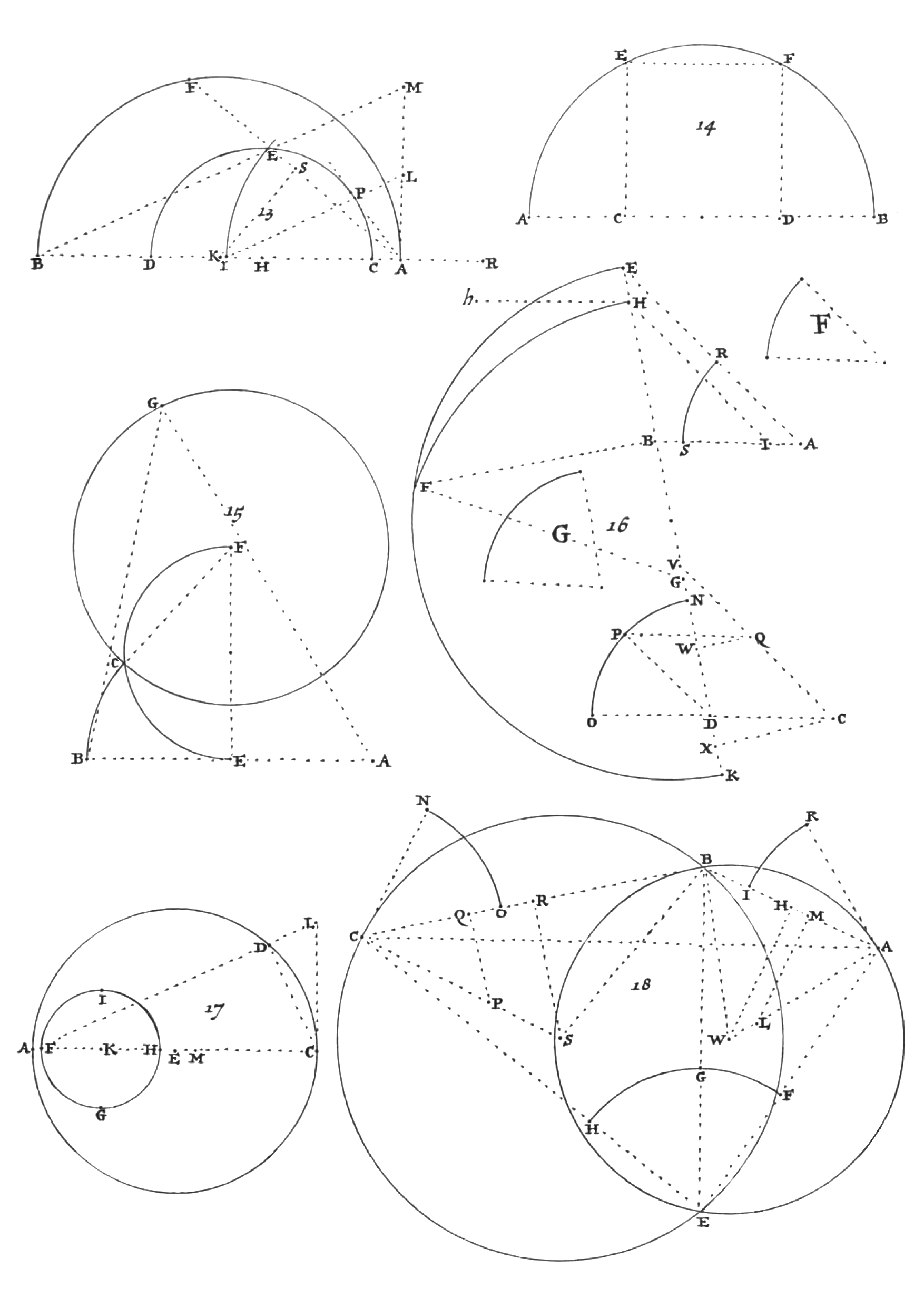}
\end{center}
\caption{Propositions 13 to 18. }\label{fig2-prop-13-18}
\end{figure}

\begin{figure}[!ht]
\begin{center}
\includegraphics[width=\textwidth]{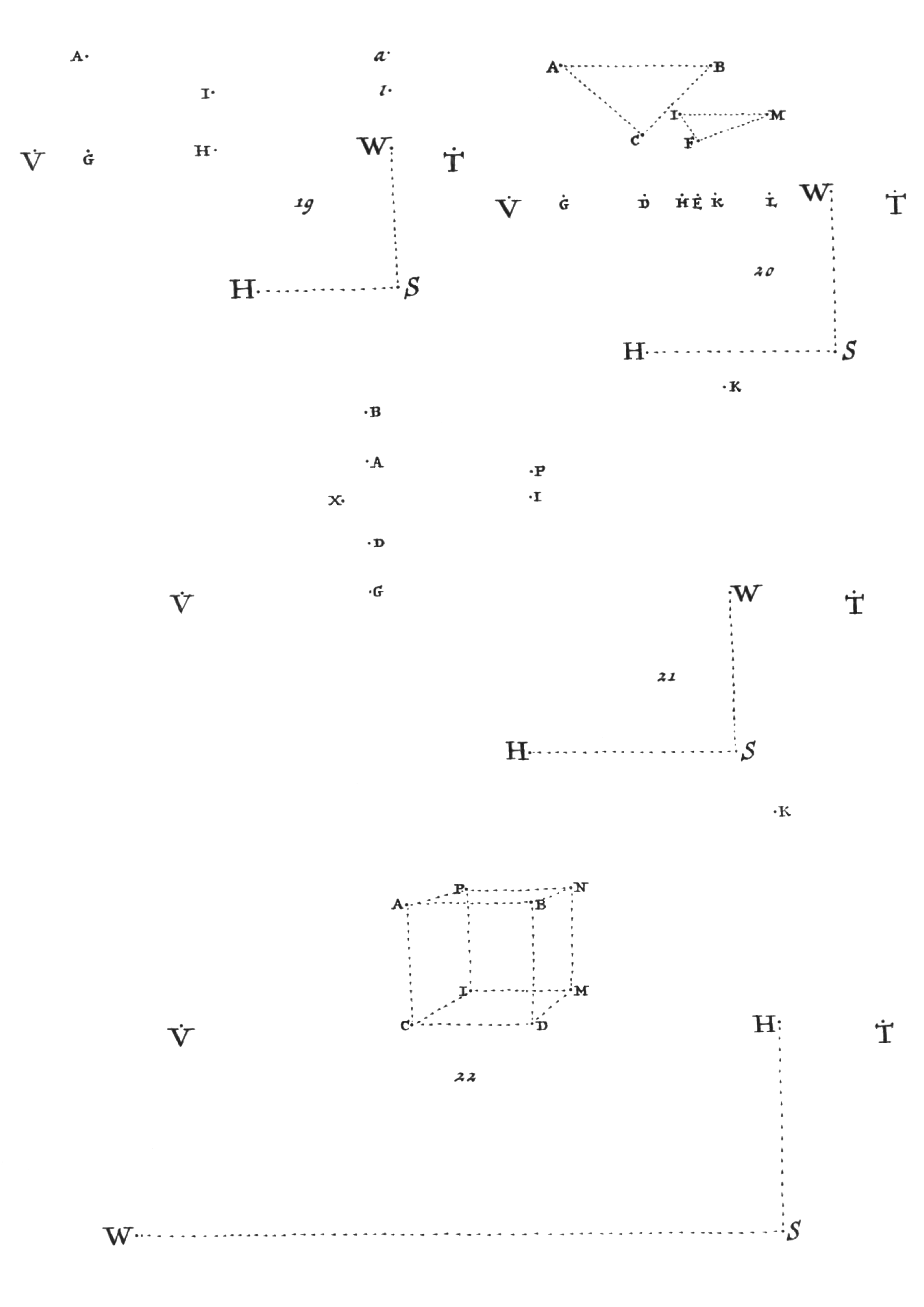}
\end{center}
\caption{Propositions 19 to 22. }\label{fig2-prop-19-22}
\end{figure}

\begin{pf}
Pull $BD$ off $BC$ from $B$, resulting in $EC$ (by 23, First Part), equal to the radius of the larger circle; then pull $EC$ from $B$ off $BD$ and\footnote{In the text it is written \emph{or}, but it clearly must be \emph{and}. See Figure \ref{fig2-prop-1-12}, page \pageref{fig2-prop-1-12}.} $BA$, resulting in $BM$ and $BN$; make a square on $BM$ (by 13, First Part), obtaining $Q$, the required centre. Now to find the radius of the smaller circle $GH$. Drop a perpendicular from $E$ to $CD$ or $CA$, obtaining $F$ (by 19, First Part); find the fourth proportional to three lines (by 31, First Part), that is, \emph{as $CE$ plus $EF$ is to $CE$, so $EF$ is to the radius $CG$, equal to $CP$}; now pull this [$CG$] from $C$ off $CD$, obtaining $GD$; make a square on $CG$, resulting in $H$, the centre; or add $CG$, equal to $EH$, to $EB$ (by 22, First Part), and the same results.
\end{pf}

\begin{proposition} 
Given a semicircle $AKC$, whose radius is $AB$ and two semicircles are described on the diameter $AC$, and the geometric mean $DK$ is also given\footnote{The diameters of the two semicircle are $CD$ and $DA$, and $DK$ is drawn perpendicular to $CA$.}. Now one wants to describe the two (largest) circles, such as $EMS$ and $HTL$, so that they should touch $DK$ (in $P$ and $O$) and also the given circles (in $M$ and $S$ and in $L$ and $T$).
\end{proposition}

\begin{pf}
Find the third proportional to two lines (by 29, First Part), that is: \emph{as $AB$ is to $AF$, so $AF$ is to $RW$}; pull this off $AF$ (by 23, First Part), then the radius of the circle $ES$, equal to the radius of the circle $HT$.
To find the centre $E$, add the radius to $AF$, obtaining $FE$ (by 22, First Part), and pull the radius off $AB$, resulting in $BE$; with these two lines make a triangle on $FB$, obtaining the centre E.
Then to find the centre $H,$ add the radius to $DG$, equal to $CG$, forming $GH$, and also pull [the radius] off $AB$, resulting in $BH$,  equal to $BE$; make a triangle on $BG$ with these two lines, then one obtains the centre $H$; now drop from $E$ and $H$ perpendiculars to $KD$ (by 19, First Part), resulting in $O$ and $P$; describe from $E$ and $H$ the two circles with the previously found radius, then they will touch, as required, in $M$, $S$ and $P$ and in $L$, $T$ and $O$.
\end{pf}

\begin{proposition} 
Given two points $A$ and $B$ above (or below) an imaginary straight line $CD$, now one wishes to describe a circle so that it should pass through the two points and touch the imaginary straight line.
\end{proposition}

\begin{pf}
Drop from $B$ a perpendicular to $CD$, falling in $E$; likewise from $A$ to $BE$, falling in $F$ (by 19, First Part); find the fourth proportional to three lines, that is \emph{as $BF$ is to $AF$, so $BE$ is to $EG$}; now add it [$EG$] to $ED$, as well as [add] $HG$, equals $GA$, to $BA$ (by 22 and 4, First Part), then find the geometric mean $GI$ between $HG$ and $GB$ (by 28, First Part), add this length [$Gl$] from $G$ to $GD$, ending in $K$; now describe a circle through the points $K$, $A$ and $B$ (by 33, First Part), then the requirement is satisfied.
\end{pf}

\begin{proposition} 
Divide a given imaginary straight line $AB$ in $C$ so that if one makes an equilateral triangle on one part $CB$, it should be as large as the square on the other part $AC$.
\end{proposition}

\begin{pf}
Make any equilateral triangle $EFG$, divide $EG$ in the middle at $H$ (by 15, First Part), add to $EG$ the geometric mean between $EH$ and $FH$ (by 28 and 22, First Part); if one now divides $AB$ in the ratio of $EG$ to its extension\footnote{That is, $EG$ plus the geometric mean.} (by 26), one obtains $C$, which is the desired point.
\end{pf}

\begin{proposition} 
One desires to find an imaginary straight line $EK$, which divides the given triangle $ABC$ into two equal parts [figures], starting from one of its sides\footnote{Lit. \emph{Een gegeven driehoeck ABC, begeertmen te deelen, als mede de drie sijden, met een ingebeelde rechte linie EK uyt een van de sijden, aen te vangen in twee gelijcke deel}. We used the German translation of this Proposition, from \cite{mohr1928}.}
\end{proposition}

\begin{pf}
Find the geometric mean between $BC$ and half of $AC$ (by 28, First Part); add the three lines $BC$, $AC$ and $AB$ together (by 22, First Part), and on their fourth part describe a semicircle (by 24 and 16, First Part); from one end point of the diameter, place the geometric mean previously found. If you then subtract and add the other side of the rectangle in the same circle from the previous fourth part, the result is $EC$ and $CK$ (by 22 and 23, First Part), as desired; so the triangle $CEK$ is equal to the figure $EBAK$ and the [sum of] lines $KC$ and $CE$ is equal to [the sum of] $EB$, $BA$ and $AK$.
\end{pf}

\begin{proposition} 
Given two semicircles,  as here $AFB$ and $CED$, whose diameters lie on the same base line, one wants  to draw the shortest imaginary straight line between their circles, which extends to one of the corners\footnote{That is, find the points $E$ and $F$ such that they are aligned with $A$ and $EF$ is the shortest possible segment.}, as here to $A$.
\end{proposition}

\begin{pf}
Add $DB$ to $CD$ from $C$ (by 22, First Part), resulting in $DR$; find the geometric mean between $AR$ and $AB$ (by 28, First Part), obtaining $AM$; with this erect a perpendicular from $A$ on $RA$ (by 21, First Part), likewise draw $AL$ [perpendicular to $AB$] equal to $AP$, which touches the circle $CPD$; find the fourth proportional to three lines, that is, \emph{as $RA$ is to $AM$, so $AL$ is to $AI$} (by 31, First Part); now add $AI$ to $RA$ from $A$, obtaining $I$; with this length $AI$ make an arc from $A$, which intersects the circle $CED$ in $E$; drop the perpendicular to $AE$ from the centre $K$ of the circle $BFA$ (by 19, First Part), falling on $S$; then make $AS$ twice as large from $A$ (by 4), it results in $F$, as desired.
\end{pf}

\begin{proposition} 
Inscribe the largest square in a semicircle.
\end{proposition}

\begin{pf}
Find the geometric mean $CD$ between the diameter $AB$ and its fifth part (by 24 and 28, First Part), which is the desired side of the square.
\end{pf}

\begin{proposition} 
From two given points $A$ and $B$ on a given imaginary straight line $AB$, one wants to draw two imaginary straight lines $AG$ and $BG$, whose squares together form with the triangle $AGB$  (construed by the given line $AB$ and the two drawn lines $AG$ and $BG$) a given ratio, as $R$ to $S$, which should be no less than 4.
\end{proposition}

\begin{pf}
Divide $AB$ in the middle at $E$ (by 15, First Part); find now the fourth proportional to three lines, that is, \emph{as $S$ is to $R$, so $AB/4$ is to $EF$}; with this [length] erect from $E$ a perpendicular to $AB$ (by 21, First Part); now describe on $EF$ a semicircle (by 16, First Part) and therein mark $EC$ equal to $EB$, and from $F$ with length $FC$ describe a circle; then if one trace from $A$ and $B$ to any point $G$ in this circle one obtains $AG$, $GB$, as desired.
\end{pf}

\begin{proposition} 
Given any two parallel imaginary straight lines, $AB$ and $CD$, to find one point $H$ outside this, from which two straight lines $HI$ and $HD$ extend at given angles $F$ and $G$ to the given straight lines $AB$ and $CD$, so that the parallelogram\footnote{Lit. \emph{rechthoeck}, ``rectangle'', which is not what Mohr intended.} $IH$, $HD$, which is enclosed by these, should be as large as a given one figure $AEBK$.
\end{proposition}

\begin{pf}
Make over $DC$ from $D$ the angles $ODN$ and $ODP$ equal to the angles $F$ and $G$, as well as on $AB$ the angle $ASR$ equal to the angle $F$ (by 8 and 23, First Part); 
from $C$ draw a parallel $CD$ to $DP$ (by 46, First Part); drop perpendiculars to $BD$ from $Q$ and $C$ (by 19, First Part) falling on $W$ and $X$; find the fourth proportional to three lines (by 31, First Part), that is, \emph{as $QW$ subtracted from $CX$ is to $XW$, so $QW$ is to $WV$}; add this to $WX$ from $W$ (by 22 First Part), arriving in $V$; now find the fourth proportional to three lines, that is, \emph{as $CD$ is to $DV$, so $AB$ is to $BE$}; add this to $BD$ from $B$, obtaining $E$; now find the geometric mean $BF$ between $EB$ and $BK$ (by 28, First Part) and halve $BD$ in $G$ (by 15, First Part)\footnote{Mohr points wrongly to Proposition 13 instead of 15.}; add  $GH$, equal to $GF$, to $GK$ from $G$; draw  from $H$ a parallel to $AB$ (by 46), say $Hh$; then if one chooses an arbitrary point on this imaginary straight line $Hh$, I choose $H$, and draw a parallel to $AE$ from it, which ends on the imaginary straight line $AB$ in $E$, that is, \emph{as $EB$ is to $AB$, so $HB$ is to $BI$} (by 31 and 23, First Part), then what is contained in $IH$, $HD$ is as large as the given figure, contained in $AE$ and $BK$, as required.
\end{pf}

\begin{proposition} 
Given a circle $ABCD$ whose radius is $AE$ and also a circle $FGHI$ whose radius is $KH$ and whose centres of gravity are the centres $E$ and $K$, find the centre of gravity of the remaining piece obtained by cutting out the smaller circle from the larger one.
\end{proposition}

\begin{pf}
Add $EC$, equal to $EA$, to $EK$ from $E$ (by 22, First Part), then mark from $C$ $CD$ equal to $HE$, which is the diameter of the smaller circle; now find the third proportional to two lines (by 30, First Part), that is, \emph{as $AD$ is to $DC$, so $DC$ is to $DL$}; extend $EK$ up to $M$ so that $KE$ has such ratio to $EM$, as $AD$ to $DL$ (which is the ratio of the remaining piece to the smaller circle), as in 31 and 22 of the First Part; this is the desired centre of gravity.
\end{pf}

\begin{proposition}\footnote{This is known as the Snellius-Pothenot Problem.} 
Given three points $A$, $B$ and $C$; whose observer stood in $E$ and found the angle $AEB$ equal to the arc $FG$ and the angle $BEC$ equal to the arc $GH$. One asks for the point $E$  where the observer was.
\end{proposition}

\begin{pf}
First find the centre $E$ from the arc\footnote{Perhaps we should call the centre by another name, $E'$ say. The arc $FGH$ need not be in its final position. The desired point $E$ appears only in the end of the construction.} $FGH$ (by 34, First Part), then set the angle $IAK$ from $A$ equal to the angle $GEF$ (by 23 and 6, First Part); set $AL$, equal to $AK$, perpendicular to $AK$ from $A$ (by 13, First Part); drop from $L$ a perpendicular to $AK$ at $M$ (by 19, First Part); find the fourth proportional to three lines (by 31, First Part), that is, \emph{as $AM$ is to $AL$, so $AH$ (equal to half of $AB$) is to $AW$, equal to $BW$}; with this length describe a circle from $W$; then set the angle $OCN$ equal to the angle $GEH$ from $C$, and set $CP$, equal to $CN$, perpendicular to $CN$ from $C$; drop from $P$ a perpendicular to $BC$ at $Q$; find the fourth proportional to three lines, that is, \emph{as $CQ$ is to $CP$, so $CR$ (equal to half of $BC$) is $CS$, equal to $BS$}; with this length describe a circle from $S$ that intersects the previous circle in $E$, which is the desired point.
\end{pf}

\begin{proposition} 
Given a point $A$ on the floor and the glass base $VT$, on which the glass\footnote{Glass is the plane of projection, floor is the plane of the figure being projected.} is perpendicular, so that both the distance and height of the observer is equal to $SH$ over $S$; now one wishes to find the picture of the point $A$.
\end{proposition}

\begin{pf}
Drop a perpendicular from $A$ and $S$ to $VT$ (by 19, First Part), at $G$ and $W$ ($SW$ is the distance the observer stands from the base of the glass); find the fourth proportional to three lines (by 31, First Part), that means: \emph{as $SW$ plus $AG$ is to $AG$, so $GW$ is to $GH$}; now put this [$GH$] with $VG$ along a straight line (by 22, First Part) and then again, \emph{as $SW$ plus $AG$ is to $AG$, so $SH$ is to $HI$}; set this perpendicular to $VH$ from $H$ (by 21, First Part), then $I$ is the desired point. But if it should happen that the point $a$ falls on the imaginary straight line $WS$, then point $a$ can be moved , parallel with $VT$, to the left or the right side of $WS$, as herein $A$, whose image is $I$; place this (as $HI$) perpendicular to $VW$ from $W$, resulting in $Wi$; then the requirement is satisfied. Or search for the fourth proportional to three lines, that is: \emph{as $aS$ is to $SH$, so $aW$ is to $Wz$}, as before.
\end{pf}

\begin{proposition} 
Given three drawn points $A$, $B$ and $C$ on the floor (or a triangle $ABC$) and the glass base $TV$, on which the glass stands vertically at right angle, in such a way that the height of the observer is $SH$; now one wishes to find their pictures.
\end{proposition}

\begin{pf}
Drop perpendiculars from $A$ and $S$ to $VT$, resulting in $G$ and $W$ and, as in the previous one, we find the image $I$ of $A$. To find the image of $B$: drop the perpendicular from $B$ to $VT$ in $K$; now look for the fourth proportional to three lines (by 31, First Part), that is, \emph{as $SW$ plus $BK$ is to $BK$, so $KW$ is to $KL$}; now set $KL$ with $VK$ from $K$ along a straight line (by 22, First Part) and then again, \emph{as $SW$ plus $BD$ in to $BK$, so $SH$ is to $LM$}; now set $LM$ from $L$ perpendicular to $VL$, then $M$ is the image of $B$. To find the image of $C$, drop the perpendicular to $VT$ from $C$, resulting in $D$; look for the fourth proportional to three lines, that is, \emph{as $SW$ plus $DC$ is to $DC$, so $DW$ is to $DE$}; now join $DE$ with $VD$ from $D$ along a straight line; then again: \emph{as $SW$ plus $DC$ is to $DC$, so $SH$ is to $EF$}; now erect $EF$ from $E$ perpendicular to $DE$; then the point $F$ is the image of point $C$ and triangle $IMF$ is the image of triangle $ABC$ as required.
\end{pf}

\begin{proposition} 
Given a drawn point $B$ above the floor, whose stand drawing is $BA$, and the glass base $VT$, on which the glass is perpendicular to the floor, also the height of the observer is $SH$; then one wishes to find the image of point $B$.
\end{proposition}

\begin{pf}
First find point $I$, which is the image of point $A$ (by 19 of this part); now erect the height $SH$ of the observer from $W$ perpendicular to $WT$ as $WK$ (by 21, First Part); drop from $A$ or $B$ a perpendicular to $VW$, resulting in $G$ (by 19, First Part); now erect $GD$, equal to $AB$,  perpendicular to $GW$  from $G$; draw a parallel to $GV$ from $I$ (by 46, First Part), as $IX$; now find the fourth proportional to three lines, that is, \emph{as $KG$ is to $GD$, so $AI$ is to $IP$}; erect $IP$ perpendicular to $IX$  from $I$; then the point $P$ is the image of point $B$.
\end{pf}

\begin{proposition} 
A cube is given, whose basis is $ABCD$; and its side $CD$ lies on a straight line with the glass base $TV$, on which the glass is perpendicular, and so that the height of the observer is $SH$; now one wants to find its picture.
\end{proposition}

\begin{pf}
First (by 20 of this part) find points $I$ and $M$, which are the images of $A$ and $B$; erect the observer's height $SH$ perpendicular to $HT$ from $H$  as in the foregoing; now find the fourth proportional to three lines (by 31, First Part), that is, \emph{as $DK$ is to $DB$, so $MK$ is to $MN$, equal to $IP$}; with then erect perpendiculars to $IM$, with $IP$, from $I$ and $M$ (by 21, First Part); then the cube $C, D, M, N, B, A, P, I$ is the desired image.
\end{pf}

\begin{proposition} 
One wants to make a horizontal sundial\footnote{This and the next proposition treat the construction of sundials with the compass alone. The reader would profit in consulting Wikipedia's article on sundials,  https://en.wikipedia.org/wiki/Sundial for detailed information on sundials.}, whose arc $BD$ is given, equal to the elevation of the Pole\footnote{In the Danish version, \emph{axel prikten}, axis point.}.
\end{proposition}

\begin{figure}[!ht]
\begin{center}
\includegraphics[width=\textwidth]{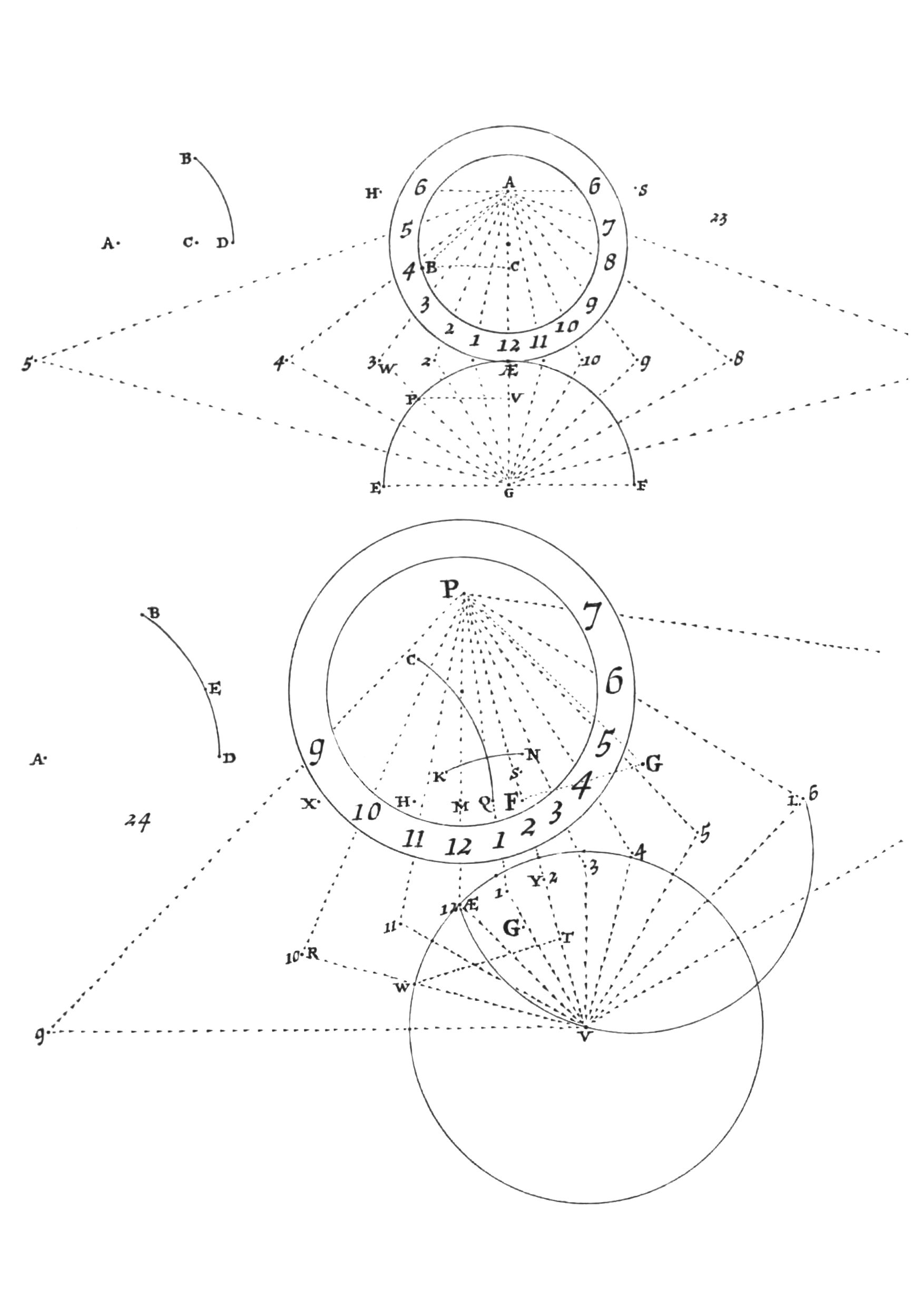}
\end{center}
\caption{Propositions 23 and 24, the Sundial }\label{fig2-prop-23-24}
\end{figure}

\begin{pf}
Find the centre [$A$] of the arc $BD$ (by 34, First Part); drop a perpendicular to $AD$  from $B$ to $C$ (by 19, First Part); now make the triangle $AC\!B$ equal to the triangle\footnote{Mohr uses the same letter for both triangles. Look at the corresponding drawing in Figure \ref{fig2-prop-23-24}, page \pageref{fig2-prop-23-24}.} $ACB$ (by 8, First Part); find the third proportional\footnote{Mohr does not make this explicit: \emph{as $AC$ is to $BC$, so $BC$ is to $C\!\mbox{\textit{\AE}}$.} He uses the letter \textit{\AE} referring to the \emph{\ae quinoctial}, or \emph{equinoctial}.} $C\!\mbox{\textit{\AE}}$ (by 30, First Part); add here $\mbox{\textit{\AE}}G$, equal to $\mbox{\textit{\AE}}B$, from $\mbox{\textit{\AE}}$ (by 22, First Part); now erect $EG$, equal to $\mbox{\textit{\AE}}G$, perpendicular to $\mbox{\textit{\AE}}G$ from $G$ (by 11, First Part) and describe a semicircle from $G$ with length $EG$ (by the Corollary of 2, First Part); divide this semicircle into twelve equal parts (by 12 and 32, First Part); drop from each of the equal parts a perpendicular to $\mbox{\textit{\AE}}G$, I take $P$, falling at $V$; now find the fourth proportional to three lines (by 31, First Part), that is, \emph{as $GP$ is to $VP$, so $G\!\mbox{\textit{\AE}}$ is to $W\!\mbox{\textit{\AE}}$}, and erect with this a perpendicular to $\mbox{\textit{\AE}}G$ at $\mbox{\textit{\AE}}$ (by 21, First Part), resulting in $W$, which is the point of the third hour and one can also find the others. Then come the hours 5, 4, 3, 2, 1, 12, 11, 10, 9, 8, 7 on the \AE quinoctial [line]. To find the sixth hour, erect from $A$ perpendiculars to $AC$ and $AG$, as  $AH$ and $AS$ (by 10 or 46, First Part) and then draw straight lines from $A$ for every hour (by the Corollary of 23, First Part); the rest shows itself and the required is met.
\end{pf}

\begin{proposition}
One wants to make a vertical sundial that deviates from South to West, which is indicated by the arc $DE$ and the increase from the Pole is equal to the arc $BD$.
\end{proposition}

\begin{pf} 
Find the centre $A$ of the arch $BD$ (by 34, First Part); now choose $FG$ at random (the central pin which should be perpendicular to the wall) and set it perpendicularly\footnote{Mohr means here \emph{perpendicular to the pin}, or in the plane of the wall. Notice that in the drawing there are two points with the same name $G$. He refers to the lowest one here.} from $F$ to $G$; make on $GF$ from $G$ the angle\footnote{Here $F$, $G$ and $N$ are collinear.} $NGK$ equal to the angle $DAE$ (by 22 and 6, First Part), which is the deviation\footnote{Deviation fro South to West.}; drop from $K$ a perpendicular to $FG$ (or $NG$) falling at $S$ (by 19, First Part); now find the fourth proportional to three lines (by 31, First Part), that is, \emph{as $SG$ is to $SK$, so $FG$ is to $FM$}; erect this perpendicularly to $FG$ from $F$; find now the third proportional of two lines  (by 30, First Part), that is, \emph{as $MF$ is to $FG$, so $FG$ is to $FL$}; add this to $FM$ from $F$, then $L$ is one of the points through which the equinoctial line passes; again add $XM$, equal to $MG$, to $MF$  from $M$ (by 22, First Part) and make on $XF$ from $X$ the sector $QXC$ equal to the sector $DAB$, which is the increase from the Pole (by 23 and 6, First Part); drop from $C$ a perpendicular to $XQ$ falling at $H$; find the fourth proportional to three lines, that is, \emph{as $XC$ is to $HC$, so $XM$ is to $MP$}; set this perpendicularly to $XM$ from $M$; now find the third proportional to two lines, that is, \emph{as $MP$ is to $XM$, so $XM$ is to $M\!\mbox{\textit{\AE}}$}; then $PM\!\mbox{\textit{\AE}}$ is the Meridian and $P$ is the point of the Pole from which all hour lines are to be drawn to the equinoctial line and {\textit{\AE}}\ is the second point of the equinoctial line; now describe a semicircle on $H\mbox{\textit{\AE}}L$ (by 16, First Part); therein mark $\mbox{\textit{\AE}}V$, equal to $\mbox{\textit{\AE}}X$, from {\textit{\AE}}, and with this length $ \mbox{\textit{\AE}}V$; describe a circle from $V$ with this length $H\!\mbox{\textit{\AE}}$, and from $\mbox{\textit{\AE}}$ start to divide the circle into 24 equal parts (by 12 and 32, First Part); now drop from $\mbox{\textit{\AE}}$ a perpendicular to $VP$, falling at $Y$, as well as from each of the same parts to $VP$, I assume from $W$, resulting in $T$ (by 19, First Part); now find the fourth proportional to three lines, that is, \emph{ as $VT$ is to $WT$, so $VY$ is to $YR$}; add this to $YL$, then $R$ is the desired point of the tenth hour; and also look for the other hours in the equinoctial line; then draw imaginary straight lines from the axis point P every hour on the equinoctial line (by the Corollary of 23, First Part), the rest is self-explanatory. If one wishes to have the heavenly symbols with it, this can also be construed.
\end{pf}

\clearpage

\section{The Figures for the Second Part}

\begin{thebibliography}{88}

\bibitem{bomhoff-dictionary}
D. Bomhoff. \textit{Nieuw Woordenboek der Nederduitsche en Engelsche Taal}, Parte II, J. F. Thieme, Nijmegen, 1851.

\bibitem{hjelmslev}
Johann Hjelmslev (ed). "Beitraege zur Lebenabschreibung von Georg Mohr (1640-1697), Kongelige Danske Videnskabernes Selskabs Skrifter, Math.-fysiske Meddelelser, 11 (1931), 3--23.

\bibitem{escofier-laouenan}
Jean-Pierre Escofier, Jean Michel Le Laou\'enan. Constructions au Compas seul. Text available at (with last access on 24th March 2020).

https://irem.univ-rennes1.fr/ouvrages-en-ligne-mathematiques-et-histoire

\bibitem{marshall-grammar}
T. Marshall. \textit{The Dutch Grammar}. Second Edition, T. Marshall \&\ Co., Rotterdam, 1854.

\bibitem{mascheroni}
Lorenzo Mascheroni. \textit{La Geometria del Compasso}. New edition, Casa Ed. Era Nova, Palermo, 1901.

\bibitem{mohr1672}
Georg Mohr. \textit{Euclides Danicus}. (Dutch version) Iacob van Helfen, Amsterdam, 1762.


\bibitem{mohr1928}
Georg Mohr. \textit{Euclides Danicus}. (Danish version) With a preface by Johannes Hjelmslev and a translation into German by P\'al. Andr. Fr. H\o st \&\ S\o n, Copenhagen, 1928.


\bibitem{schouten}
Franciscus van Schouten. \textit{Tractaet der Perspective, ofte Schynbaere Teycken-Konst}. Gerrit van Goedesbergh, Amsterdam, 1660. Available in archive.org.

\bibitem{strand-tortzen}
Hanne Eggert Strand, Chr. Gorm Tortzen. Euclides Danicus -- en studie i Georg Mohrs honnette ambition. Festskrift til Bent Christensen i anledning af hans 70 {\aa}rs f{\o}dselsdag. AIGIS, Vol. 16, No. 1, 2016. Supplementum V. Available in
http://aigis.igl.ku.dk/aigis/BC70/Euclid.pdf (last access, April, 10th, 2018).

\end{thebibliography}
\end{document}